\documentclass{amsart}

\usepackage{amsmath,amssymb,yhmath}

\newtheorem{lemma}{Lemma}
\newtheorem{prop}[lemma]{Proposition}
\newtheorem{cor}[lemma]{Corollary}
\newtheorem{thm}[lemma]{Theorem}

\newtheorem{thm?}[lemma]{Theorem?}

\newtheorem{conj}{Conjecture}

\newtheorem*{unthm}{Theorem}
\newtheorem{fact}{Fact}
\newtheorem{prob}{Problem}

\begin{document}
\title{Torsion Points on Elliptic Curves with Complex Multiplication}
\author{Pete L. Clark}
\author{Brian Cook}
\author{James Stankewicz}

%\author{PC}
%\author{BC}
%\author{AR}
%\author{JS}
%\author{NW}
%\author{SW}
%\author{BW}

\address{Department of Mathematics \\ Boyd Graduate Studies Research Center \\ University 
of Georgia \\ Athens, GA 30602-7403 \\ USA}
\email{pete@math.uga.edu}
\email{bcook@math.uga.edu}
\email{stankewicz@gmail.com}

\maketitle

\newcommand{\F}{\mathbb{F}}
\newcommand{\et}{\textrm{\'et}}
\newcommand{\ra}{\ensuremath{\rightarrow}}
\newcommand{\FF}{\F}
\newcommand{\Z}{\mathbb{Z}}
\newcommand{\N}{\mathcal{N}}
\newcommand{\ch}{}
\newcommand{\R}{\mathbb{R}}
\newcommand{\PP}{\mathbb{P}}
\newcommand{\C}{\mathbb{C}}
\newcommand{\Q}{\mathbb{Q}}
\newcommand{\tpqr}{\widetilde{\triangle(p,q,r)}}
\newcommand{\ab}{\operatorname{ab}}
\newcommand{\Aut}{\operatorname{Aut}}
\newcommand{\gk}{\mathfrak{g}_K}
\newcommand{\gq}{\mathfrak{g}_{\Q}}
\newcommand{\OQ}{\overline{\Q}}
\newcommand{\Out}{\operatorname{Out}}
\newcommand{\End}{\operatorname{End}}
\newcommand{\Gon}{\operatorname{Gon}}
\newcommand{\Gal}{\operatorname{Gal}}
\newcommand{\CT}{(\mathcal{C},\mathcal{T})}
\newcommand{\ttop}{\operatorname{top}}
\newcommand{\lcm}{\operatorname{lcm}}
\newcommand{\Div}{\operatorname{Div}}
\newcommand{\OO}{\mathcal{O}}
\newcommand{\rank}{\operatorname{rank}}
\newcommand{\tors}{\operatorname{tors}}
\newcommand{\IM}{\operatorname{IM}}
\newcommand{\CM}{\operatorname{CM}}
\newcommand{\Frac}{\operatorname{Frac}}
\newcommand{\Pic}{\operatorname{Pic}}
\newcommand{\coker}{\operatorname{coker}}
\newcommand{\Cl}{\operatorname{Cl}}
\newcommand{\loc}{\operatorname{loc}}
\newcommand{\GL}{\operatorname{GL}}
\newcommand{\PSL}{\operatorname{PSL}}
\newcommand{\Frob}{\operatorname{Frob}}
\newcommand{\Hom}{\operatorname{Hom}}
\newcommand{\Coker}{\operatorname{\coker}}
\newcommand{\Ker}{\ker}

\noindent

\begin{abstract}
We present seven theorems on the structure of prime order torsion points on CM elliptic curves defined over number 
fields.  The first three results refine bounds of Silverberg and Prasad-Yogananda by taking into account the class number of the CM order and the splitting of the prime in the CM field.  In many cases we can show that our refined bounds are 
optimal or asymptotically optimal.  We also derive asymptotic upper and lower bounds on the least degree of a CM-point on 
$X_1(N)$.  Upon comparison to bounds for the least degree for which there exist infinitely many rational points 
on $X_1(N)$, we deduce that, for sufficiently large $N$, $X_1(N)$ will have a rational CM point of degree smaller than 
the degrees of at least all but finitely many non-CM points.
\end{abstract}

\section{Introduction}

\subsection{Notation} \textbf{} \\ \\ \noindent
For $d \in \Z^+$, we define the following quantities: \\ \\
$T(d)$: the supremum of the orders of the groups $E(K)[\tors]$ as $K$ ranges over all 
number fields of degree $d$ and $E$ ranges over all elliptic curves defined over $K$. 
\\ \\
$N(d)$: the supremum of all orders of $K$-rational torsion points $P \in E(K)$, with $K$ and 
$E$ varying as above.
\\ \\
$P(d)$: the supremum of all prime orders of $K$-rational torsion points $P \in E(K)$, with $K$ and 
$E$ varying as above.
\\ \\
We shall have occasion to consider analogues $T_*(d)$, $N_*(d)$, $P_*(d)$ of the above quantities, which are 
defined by restricting to some subset of elliptic curves $E_{/K}$.  Specifically we will be interested 
in the set of all elliptic curves with integral modulus $j(E)$ and also the set of all elliptic 
curves with complex multiplication.    
%\\ \\  
%For brevity, we use the abbreviation 
%\[\Z(a_1,\ldots,a_n) := \Z/a_1\Z \times \ldots \times 
%\Z/a_n\Z. \]
\subsection{Background on torsion} \textbf{} \\ \\ \noindent
Since the torsion subgroup of an elliptic curve over a number field is a finite abelian group with at most 
two generators, we have 
\begin{equation}
\label{FIRSTINEQ}
P(d) \leq N(d) \leq T(d) \leq N(d)^2. 
\end{equation}
\noindent
The \textbf{uniform boundedness theorem} of L. Merel \cite{Merel} asserts $T(d) < \infty$ for all 
$d \in \Z^+$.  Using (\ref{FIRSTINEQ}), the finiteness of $P(d)$ and $N(d)$ follows immediately.  
\\ \\
Merel's proof gives an explicit upper bound on $T(d)$, which was then improved by work of Merel, Oesterl\'e and Parent.  For instance, Parent showed \cite{Parent99} that 
if a power $p^a$ of a prime $p > 3$ divides the order of the torsion subgroup of an elliptic curve over a degree 
$d$ number field, then \[p^a \leq 65(3^d-1)(2d)^6. \]  However, it is a ``folk conjecture'' that there exists a 
constant $\alpha$ such that $T(d) = O(d^{\alpha})$: thus it seems that Merel's bounds are a full exponential 
away from the truth.  In fact, we record here a more precise conjecture:
\begin{conj}
\label{GUESSCONJ} \textbf{} \\
There is a $C_2 > 0$ such that $T(d) \leq C_2 d \log \log d$ for all $d \in \Z^+$. \\
\end{conj}
\noindent
Conjecture \ref{GUESSCONJ} is very close to being the most ambitious conceivable one: we shall show 
(Theorem 6) that there is a positive constant $C_1$ and a strictly increasing sequence $\{d_n\}_{n=1}^{\infty}$ 
of positive integers such that $T(d_n) > C_1 d_n \sqrt{\log \log d_n}$ for all $n$.  
\\ \\
Unfortunately it is not currently tenable to seek numerical confirmation for
Conjecture \ref{GUESSCONJ}a).  The only values of $d$ for which any of $T(d)$, $N(d)$, $P(d)$ are known are: 
\\ \\
$T(1) = 16$, $N(1) = 12$, $P(1) = 7$ (\cite{Mazur}). \\ \\
$T(2) = 24$, $N(2) = 18$, $P(2) = 13$ (\cite{Kamienny86}, \cite{Kamienny92},
\cite{KenkuMomose}).  
\\ \\
$P(3) = 13$ (\cite{Parent}).
\\ \\
Since further direct computation of these quantities is out of current reach, it seems that one 
must find some more tractable sub-problem and examine the extent to which it is 
representative of the general case.  \\ \indent One approach is to concentrate on the case of elliptic 
curves with algebraic integral $j$-invariant (henceforth \textbf{integral modulus}).  In this case we 
write $T_{\IM}(d), \ N_{\IM}(d)$, $P_{\IM}(d)$ for the order, exponent and largest prime dividing the order 
of an elliptic curve $E$ with integral modulus defined over any number field of degree $d$.  For such curves the uniform boundedness is much 
easier to prove.  Moreover, in the integral modulus case the computation of all possible torsion subgroups over 
$\Q$ was done by G. Frey in 1977 \cite{Frey}.  Analogous computations in higher degree are 
significantly more difficult and have been the subject of several papers of H. Zimmer and his collaborators: the 1976 paper \cite{Zimmer0} lays foundations by giving a generalization of the Lutz-Nagell 
restrictions on torsion points to arbitrary number fields; the 1989 paper \cite{Zimmer1} enumerates 
the torsion subgroups of elliptic curves with integral modulus over quadratic fields ($d=2$); special kinds of cubic fields ($d = 3$) 
were considered in 1990 \cite{Zimmer2} and the case of a general 
cubic field was completed in 1997 \cite{Zimmer3}; only a very restricted class of quartic fields has 
ever been considered, so already the case $d = 4$ seems to be out of reach.
\\ \indent However, Hindry and Silverman have shown \cite{HS} that 
\begin{equation}
\label{HSONE}
\forall d \in \Z^+, T_{\IM}(d) \leq 1977408 d \log d, 
\end{equation}
\begin{equation}
\label{HSTWO}
\forall d \geq 25, T_{\IM}(d) \leq 498240 d \log d. 
\end{equation}
\indent
Another idea is to search for all finite groups which arise as the torsion subgroup 
of \emph{infinitely} many elliptic curves defined over number fields of degree $d$. In this case the computations in degree up to $d = 4$ have been 
done by Jeon, Kim, Park and Schweizer \cite{JKS}, \cite{JK}, \cite{JK06}, and reasonably good asymptotic bounds can be obtained 
by applying theorems of Faltings and Abramovich.  This work is described in some detail below.    
\\ \\
In this paper we shall usually restrict to elliptic curves with complex multiplication.  This is a very special subclass of the class of integral moduli curves, comprising 
for each degree $d$ only finitely many j-invariants (but infinitely many nonisomorphic \emph{twists} for a given $j$-invariant).  
Accordingly, we are able to derive more precise results than in the general case.  We also take 
up the task of relating the special case of CM points to the general case -- not definitively, 
of course, but in a depth and level of detail which we feel deserves a place in the literature on the subject.  

\subsection{Prior results} \textbf{} \\ \\ 
\noindent
Let $F$ be a field of characteristic $0$ and $E_{/F}$ an elliptic curve.  We say that $E$ has \textbf{complex multiplication} 
(henceforth \textbf{CM}) if the ring $\End E$ of endomorphisms of $E$ defined over an algebraic closure $\overline{F}$ of $F$ is strictly larger 
than $\Z$.  In this case, $\End^0(E) := \End(E) \otimes_{\Z} \Q$ is an imaginary quadratic field 
$\Q(\sqrt{D})$ and $\End(E)$ is an \textbf{order} in $\End^0(E)$.  
%(There are two isomorphisms from $\End^0(E)$ 
%to $\Q(\sqrt{D})$, but the complex analytic theory provides a distinguished one, via the representation of 
%$\End^0(E)$ on the space of holomorphic differentials of $E$: c.f. \cite[Prop. II.1.1]{SilvermanII}.)  Conversely, if $\mathcal{O}$ is 
%any order in an imaginary quadratic field $K$, we will say that an elliptic curve $E$ has $\mathcal{O}$-CM if 
%$\End(E) = \mathcal{O}$.
\\ \indent As alluded to above, we write $T_{\CM}(d)$, $N_{\CM}(d)$, $P_{\CM}(d)$ for, respectively, 
the largest order, exponent and prime dividing the order of any CM elliptic curve defined over any number field 
of degree $d$.
\\ \indent
The $j$-invariant of a CM elliptic curve is an algebraic integer \cite[Thm. II.6.1]{SilvermanII}, so that 
(\ref{HSONE}), (\ref{HSTWO}) we have $\# E(F)[\tors] = O(d \log d)$.  If we restrict to the order of a single 
torsion point -- i.e., to $N_{\CM}(d)$ rather than $T_{\CM}(d)$ -- we can do qualitatively better: one knows that 
$N_{\CM}(d) = o(d \log d)$.  More precisely:
\begin{unthm}(Silverberg \cite{Silverberg}, Prasad-Yogananda \cite{PY}) 
\label{SPYTHM} Let $F$ be a number field of degree $d$, and 
let $E_{/F}$ be an elliptic curve with complex multiplication by an order 
$\mathcal{O}$ in the imaginary quadratic field $K$.  Let $ w = w(\mathcal{O}) = 
\# \mathcal{O}^{\times}$ (so $w = 2$, $4$ or $6$) and let $e$ be the maximal order of an element of
$E(F)[\tors]$.  Then: \\
a) $\varphi(e) \leq w d$ ($\varphi$ is Euler's totient function). \\
b) If $F \supseteq K$, then $\varphi(e) \leq \frac{w}{2} d$. \\
c) If $F$ \textbf{does not} contain $K$, then $\varphi(\# E(F)[\tors]) \leq w d$.  
\end{unthm}
\noindent
Applying the theorem necessitates separate consideration of three cases:
\\ \\
Case 1: $\mathcal{O} = \Z[\frac{1+\sqrt{-3}}{2}]$, of discriminant $-3$, 
which has $w(\mathcal{O}) = 6$.  We get 
\begin{equation}
\label{SPY1}
\varphi(e) \leq 6d.
\end{equation}
Case 2: $\mathcal{O} = \Z[\sqrt{-1}]$, of discriminant $-4$, which 
has $w(\mathcal{O}) = 4$.  We get 
\begin{equation}
\label{SPY2}
\varphi(e) \leq 4d. 
\end{equation}
Case 3: For every other order we have 
$w(\mathcal{O}) = 2$.  We get 
\begin{equation}
\label{SPY3} 
\varphi(e) \leq 2d.
\end{equation}
Let us call (\ref{SPY1}), (\ref{SPY2}) and (\ref{SPY3}) the \textbf{SPY bounds}.  
\\ \\
Recall the classical result
$\varphi(N) \gg \frac{N}{\log \log N}$ (e.g. \cite[Thm. 328]{HardyWright}).  From this and the SPY bounds we deduce that 
there exists a constant $C$ such that
\begin{equation}
\label{LOGLOGEQ}
 N_{\CM}(d) \leq C d \log \log d.
\end{equation}
This improves upon what one gets by applying (\ref{HSONE}):
\[N_{\CM}(d) \leq N_{\IM}(d) \leq T_{\IM}(d) \leq 1977408 d \log d. \]
Theorem \ref{THEOREM7} below asserts $N_{\CM}(d) \neq o(d \sqrt{\log \log d})$, so that our understanding of the true lower 
order of magnitude of $N_{\CM}(d)$ is rather good.  On the other hand, it is 
vexing that we cannot get any improvement on 
\[T_{\CM}(d) \leq T_{\IM}(d) = O(d \log d)\] by applying the methods of SPY, or indeed 
by any other means that we know.

\subsection{Computational results} \textbf{} \\ \\ 
\noindent
We briefly report on some calculations done by the University of Georgia Number Theory VIGRE Research Group, which 
has implemented an algorithm (c.f. \cite{Clark1}) to do the following:
given a positive integer $d$, compute the complete list of isomorphism classes 
of finite abelian groups which arise as the full torsion subgroup of some CM elliptic curve with 
defined over any number field of degree $d$.  \\ \indent 
This algorithm requires knowledge 
of the CM j-invariants (more precisely, their minimal polynomials) of degree $d'$ 
\emph{strictly dividing} $d$, so in full generality requires an enumeration of the set of imaginary 
quadratic fields with any given class number, i.e., an effective solution of the \textbf{Gauss 
class number problem}.  Work of Watkins \cite{Watkins} gives a solution to this problem up to class number $100$,
so the data from \emph{ibid.} enable us, in theory, to run the algorithm for all degrees up to 
$d = 201$.  But in fact this is much more class number data than we have been able to use: one of the steps in our algorithm is the computation 
of an explicit polynomial $P_N(x,y) = 0$ which (birationally) defines the modular curve $X_1(N)$, 
a computation which became prohibitively expensive for us around $N = 79$.  
The complete list of possible torsion subgroups of CM elliptic curves defined over any degree $d$ number field has 
been computed by our VIGRE research group for $1 \leq d \leq 13$ (but will be described elsewhere).  The case of $d = 1$ 
is a 1974 result of L. Olson \cite{Olson}.  For $d = 2$ and $3$ the results are subsumed by the calculations of \cite{Zimmer1}, 
\cite{Zimmer3}.  To the best of our knowledge the cases $4 \leq d \leq 13$ had not been computed before.
\\ \\
Upon restriction from $T_{\CM}(d)$ to $P_{\CM}(d)$, the above problem can be rephrased as follows: for a 
fixed $d$, find all prime numbers $N$ such that the modular curve $X_1(N)$ has a CM point of degree $d$.  It is natural to consider also the 
following ``converse problem'': for fixed prime $N$, find the smallest degree of a CM point on $X_1(N)$.  Our 
algorithm works equally well on this converse problem, and we present the solution, for all $N \leq 79$, 
in the following table:\footnote{Some preliminary calculations were done by the first author.  The calculations were rechecked 
and completed by Steve Lane, who also pointed out -- several times -- an error in the preliminary calculations at $N = 11$, 
which turned out to be very interesting and significant.}  
\\ \\
\textbf{TABLE 1}
\\ 
$N = 2$: $d = 1$, $D = -3, -4,-7, -8, -12, -16, -28$ \\ 
$N = 3$: $d = 1$, $D = -3, -12, -27$ \\
%$N = 4$: $d = 1$, $D = -1, -4$ \\ 
$N = 5$: $d = 2$, $D = -4$\\  
%$N = 6$: $d= 1$, $D= -3, -12$ \\ 
$N = 7$: $d= 2$, $D = -3$ \\
$N = 11$: $d = 5$, $D = -11$. \\
$N = 13$: $d = 4$, $D = -3$. \\
$N = 17$: $d = 8$, $D = -4$. \\
$N = 19$: $d = 6$, $D = -3$.  \\
$N = 23$: $d = 22$, $D = -7, -11, -19, -28, -43, -67$. \\
$N = 29$: $d = 14$, $D = -4$. \\
$N = 31$: $d = 10$, $D = -3$. \\
$N = 37$: $d = 12$, $D = -3$. \\
$N = 41$: $d = 20$, $D = -4$. \\
$N  = 43$: $d = 14$, $D = -3$. \\
$N = 47$: $d = 46$, $D = -11, -19, -43, -67, -163$. \\ 
$N = 53$: $d = 26$, $D = -4$. \\
$N = 59$: $d = 58$, $D = -8, -11, -43, -67$. \\
$N = 61$: $d = 20$, $D = -3$. \\
$N = 67$: $d = 22$, $D = -3$. \\
$N = 71$: $d = 70$, $D = -7, -11, -28, -67, -163$. \\
$N = 73$: $d = 24$, $D = -3$. \\
$N = 79$: $d = 26$, $D = -3$.
\\ \\
\noindent
%In particular, 
%among entries in the table the SPY bounds are attained (i.e., we have \emph{equality}) in precisely 
%the following cases:
%\\ \\
%$\bullet$ $d = 1$, $D = -27$, $\varphi(e(T)) = \varphi(3) = 2 \cdot d$. \\
%$\bullet$ $d = 1$, $D = -12$, $\varphi(e(T)) = \varphi(6) = 2 \cdot d$. \\
%$\bullet$ $d = 5$, $D = -43$, $\varphi(e(T)) = \varphi(11) = 2 \cdot d$. \\ \\
Looking through the data one observes that most, but not all, of the time, the SPY bounds are not sharp, so it 
is natural to ask for refinements.  In the next section we shall present several such results.  Theorem \ref{NEWTHM} refines the SPY 
bounds, by including a factor of the class number $h(D)$ as well as giving a much larger lower bound in case 
$(\frac{D}{N}) = -1$.  Theorem \ref{THEOREM4} gives conditions under which one gets an extra factor of $2$ in 
the SPY-type bounds.  Moreover, for $N$ sufficiently large compared to $D$, the bounds of Theorem \ref{THEOREM4} 
are optimal.

\subsection{Theoretical results I: Optimal bounds on prime order torsion points} 
\noindent

\begin{thm}
\label{JZEROTHM} \textbf{} \\
a) For every prime $N \equiv 1 \pmod 3$, there exists an elliptic curve $E$ over a 
number field $K$ of degree $\frac{N-1}{3}$, with $j(E) = 0$, and with a $K$-rational $N$-torsion 
point. \\
b) There exists an absolute constant $N_0$ such that for all primes $N \geq N_0$: \\
(i) if $X_1(N)$ has a CM point of degree $d$, then $d \geq \frac{N-1}{3}$; \\
(ii) if $X_1(N)$ has a CM point of degree $d < \frac{N-1}{2}$ then $d = \frac{N-1}{3}$ and $j(E) = 0$.
\end{thm}
\noindent
Remark 1.1: The data suggests that it may be possible to take $N_0 = 5$.

\begin{thm}
\label{NEWTHM}
Let $\mathcal{O}_K$ be the maximal order in $K = \Q(\sqrt{D})$, $F$ a number field, and $E_{/F}$ an elliptic curve 
with $\mathcal{O}_K$ multiplication.  Let $w(K) = \# \mathcal{O}_K^{\times}$.  Suppose that $E(F)[\tors]$ contains an element of odd prime order $N$.  
Define $\delta(F,K)$ to be $1$ if $K$ is contained in $F$ and $2$ otherwise. \\
a) $(\frac{D}{N}) = 1$, then 
\[(N-1) \cdot \frac{\delta(F,K) h(K)}{w(K)} \ | \ [F:\Q]. \]
b) If $(\frac{D}{N}) = 0$, then 
\[ (N-1) \cdot \frac{(3-\delta(F,K))h(K)}{w(K)} \ | \ [F:\Q]. \]
c) If $(\frac{D}{N}) = -1$, then 
\[(N^2-1) \cdot \frac{h(K)}{w(K)} \ | \ [F:\Q]. \]
\end{thm}
\noindent
It is interesting to compare this with the SPY-bounds.  Our Theorem \ref{NEWTHM} is more special in 
that it only applies to the case of torsion points of odd prime order (although we believe the methods should 
generalize to arbitrary $N$).  In the case of prime $N$, it does not \emph{strengthen} the SPY-bound -- 
indeed, both bounds agree in the case when $N \ | \ D$, but it significantly \emph{refines} the SPY-bounds, 
making clear that they are in some sense a ``worst case scenario.''  
%On the other hand, we will later show that 
%in all cases we can have equality in these bounds, which shows in particular that if one desires a bound on the degree of $[F:\Q]$ depending 
%only on $N$, the SPY-bounds are best possible.
%
\begin{thm}
\label{THEOREM4}
Let $\mathcal{O}$ be an order in the field $K = \Q(\sqrt{D})$, $w(\mathcal{O})$ be the cardinality 
of its unit group and $h(\mathcal{O}) = \# \Pic(\mathcal{O})$ its class number.  Then:\\
a) For every odd prime $N$ which splits in $K$, there exists an $\mathcal{O}$-CM elliptic curve defined over a number field of degree 
$2(N-1) \cdot \frac{h(\mathcal{O})}{w(\mathcal{O})}$ with a rational $N$-torsion point. \\
b) There is an $N_0 = N_0(D)$ such that for $N \geq N_0$, the least degree of 
an $\mathcal{O}(D)$-CM point on $X_1(N)$ is $2(N-1) \cdot \frac{h(\mathcal{O})}{w(\mathcal{O})}$ if $N$ splits in $K$ and 
$\left( {N^2-1}\right) \frac{h(\mathcal{O})}{w(\mathcal{O})}$ otherwise.
\end{thm}
\noindent
Remark 1.2: Taking $\mathcal{O}$ to be the quadratic order of discriminant $-3$ in Theorem \ref{THEOREM4}a), we 
recover Theorem \ref{JZEROTHM}a).  The other parts of Theorem \ref{JZEROTHM} are quick consequences of 
Theorem \ref{THEOREM4} together with the SPY-bounds, but it seems worthwhile to call attention to the extremal behavior coming 
from the quadratic orders with nontrivial units.

\subsection{Theoretical results II: CM points of small degree on $X_1(N)$} \textbf{} \\ \\
Throughout this section $N$ denotes a prime number different from $2$ and $3$.
\\ \\
Define $d_{\CM}(N)$ to be the least degree of a CM point on $X_1(N)$.
\\ \\
Theorem \ref{JZEROTHM} shows that the smallest (resp. second smallest) possible degree of a CM point on $X_1(N)$ is 
$\frac{N-1}{3}$ (resp. $\frac{N-1}{2}$), and shows that this degree can be attained iff $N \equiv 1 \pmod 3$ (resp. 
$N \equiv 1 \pmod 4$).  In particular, as $N$ ranges over all primes $N$ which are \textbf{not} $11 \pmod{12}$, the least 
degree of a CM point on $X_1(N)$ is linear in $N$.  Notice that the excluded set of primes $N \equiv 11 \pmod{12}$ has 
density $\frac{1}{4}$ in the set of all primes.  By Theorem \ref{NEWTHM}, the problem of bounding the 
upper order of $d_{\CM}(N)$ as $N$ ranges over prime numbers, comes down to finding, for a given prime $N$, an imaginary quadratic field $\Q(\sqrt{D})$ such that $(\frac{D}{N}) = -1$ and with class 
number $h(D)$ as small as possible.  By applying what is known about 
these elementary -- but difficult! -- analytic problems, we arrive at the following result.
\begin{thm}
\label{LEASTDEGTHM}
a) For any $\epsilon > 0$, there exists $C_{\epsilon}$ such that for any prime $N$, 
the curve $X_1(N)$ has a CM point of degree at most $C_{\epsilon} N^{1+c/2+\epsilon}$, where 
$c/2 = \frac{1}{8} e^{\frac{-1}{2}} \approx .078$. \\
b) Assuming the Generalized Riemann Hypothesis (GRH), the least degree of a CM point on 
$X_1(N)$ is $O( N \log N \log \log N)$. 
\end{thm}
\noindent
However, $d_{\CM}(N)$ is not bounded by a linear function of $N$.
\begin{thm}
\label{NONLINEARTHM}
For any $C > 0$, there is a positive density set $\mathcal{P}$ of prime numbers such that for all 
$N \in \mathcal{P}$, the least degree of a CM point on $X_1(N)$ exceeds $CN$.  
\end{thm}
\noindent
%Remark: The restriction to admissible points is purely a technical one; we expect that part c) is true 
%for all CM points.

%
\begin{thm}
\label{THEOREM7}
a) There exists $C > 0$ such that for any $F/\Q$ with $[F:\Q] = d$ and any
CM elliptic curve $E_{/F}$, one has  $\exp(E(F)[\tors]) \leq C d \log \log d$. \\
b) There exists a sequence $F_n$ of number fields, of degree $d_n = [F_n:\Q]$ tending to infinity, 
and CM elliptic curves $E_n/F_n$ such that 
\[ \exp (E_n(F_n)[\tors]) \gg d_n \sqrt{\log \log d_n}. \]
\end{thm}
\noindent
We have already seen that part a) is a consequence of the SPY bounds; we repeat it here for the sake of 
parallelism.  Neither is part b) very difficult: all in all Theorems \ref{LEASTDEGTHM} and \ref{NONLINEARTHM} seem 
to lie significantly deeper.  

\subsection{Theoretical results III: small degree points on $X_1(N)$: comparison with non-CM case} \textbf{} \\ \\
The overarching problem is to understand all points of degree $d$ on the family 
of modular curves $X_1(N)$.  Merel's theorem asserts that for fixed $d$ the set of all such points on $X_1(N)$ is finite, 
so it is natural to enumerate this list.  Conversely, one can fix $N$ and ask for the least degree of a noncuspidal point 
on $X_1(N)$.  In the previous section we presented results giving rather tight estimates on the least degree of a 
noncuspidal CM point.  Therefore the key issue is: how many \textbf{non}-CM points are there of small degree?  
\\ \\
The next result gives a precise sense in which $d \approx N^2$ is 
the threshold between small degree and large degree:
\begin{thm}
\label{THEOREM5} \textbf{} Let $N > 3$ be a prime number.  Then: \\
a) The set of points of $X_1(N)$ of degree less than $\lceil \frac{7}{3200} (N^2-1) \rceil$ 
is finite.  Assuming Selberg's eigenvalue conjecture the bound can be improved to 
$\lceil \frac{1}{384} (N^2-1) \rceil$.  \\
b) The set of points of $X_1(N)$ of degree at most $\frac{N^2-12N+11}{12}$ is infinite.
\end{thm}
\noindent
Remark 1.3: The proof of part a) uses deep theorems of Faltings, Frey and Abramovich, but the deduction itself is now routine.  
Essentially the same result appears as \cite[Cor. 1.4]{JKS}, the only difference being that we get a sharper bound 
by restricing to prime $N$.  Part b) is much more elementary.  Nevertheless, it is in the spirit of this paper to 
pursue quantitative rather than just qualitative results, and in this regard the fact that we can compute the ``threshold'' 
value of $d$ sharply to within a factor of $32$ seems interesting.  For instance, it raises the question of whether the 
truth lies closer to $\frac{1}{384}N^2$ or to $\frac{1}{12}N^2$.  
\\ \\
Remark 1.4: Selberg's eigenvalue conjecture states that for a modular curve $Y(\Gamma) := \Gamma \backslash \mathcal{H}$ associated 
to a congruence subgroup $\Gamma \subset PSL_2(\Z)$, the least positive eigenvalue $\lambda_1$ of the hyperbolic Laplacian 
on $Y(\Gamma)$ satisfies $\lambda_1 \geq \frac{1}{4}$.  Selberg himself showed $\lambda_1 \geq \frac{3}{16}$; in 1994, Luo, 
Rudnick and Sarnak showed $\lambda_1 \geq \frac{21}{100}$; this the bound we use in our unconditional estimate.  As of 
this writing, the best known estimate on $\lambda_1$ is due to Kim and Sarnark: $\lambda_1 \geq \frac{975}{4096} > 0.238$.  
Thus the improvement in the upper bound of part a) gained by assuming Selberg's conjecture is small compared to the 
discrepancy between the upper bound of part a) and the lower bound of part b), so ought not to be the focus of our concern.
%\\ \\
%By Theorem \ref{THEOREM5} there are only finitely many points on $X_1(N)$ of 
%degree at most $C_1 N^2$, and by Theorem \ref{LEASTDEGTHM} there is a CM point on $X_1(N)$ of degree at most 
%$C_2 N^{1+\frac{c}{2}+.5} = C_2 N^{1.576}$.  It follows that for all sufficiently large $N$ there are at most finitely 
%many non-CM points on $X_1(N)$ of smaller degree than any CM point.  
\\ \\
Application: For $N = 127$ the least degree of a rational CM point is $42$, whereas -- assuming 
Selberg's eigenvalue conjecture --  the bound of Theorem \ref{THEOREM5}a) gives that there are only finitely many points (if any, of course!) on $Y_1(127)$ 
of any smaller degree.  For all larger $N \equiv 1 \pmod 3$, the set of points whose degree is 
less than or equal to the minimal degree of a CM point is finite.  
\\ \indent
On the other hand, Theorem \ref{THEOREM5}b) guarantees that there are infinitely many points of degree less than the 
smallest CM point for $N \leq 13$.  When $N = 17$ the bound ensures infinitely many points of 
degree at most $8$, and the table above shows that the least degree of a rational CM point is $8$.  But in fact there 
exists a degree $4$ map from $X_1(17)$ to the projective line, so that there are infinitely many rational points of degree at most $4$.  This suggests that there is 
room for improvement in the bound of Theorem \ref{THEOREM5}b).
\\ \\
Write $d_{\CM}(N)$ for the least degree of a CM point on $X_1(N)$ and $d_{\infty}(N)$ for the least degree $d$ such 
that $X_1(N)$ has infinitely many points of degree at most $d$.  Then by Theorem \ref{LEASTDEGTHM}, $d_{\CM}(N) = O(N^{1.078\ldots})$ 
whereas $d_{\infty}(N) \geq \lceil \frac{7}{3200} (N^2-1) \rceil - 1$.  It follows that there exists a prime $N_0$ 
such that $d_{\CM}(N_0) \geq d_{\infty}(N_0)$ and $d_{\CM}(N) < d_{\infty}(N)$ for all $N > N_0$.  In other words, 
for all sufficiently large primes, there are only finitely many points on $X_1(N)$ of degree smaller than that of any 
CM point.  \\ \\
The prime $N_0$ of the previous paragraph is effectively computable.  Indeed, B. Cook and A. Rice are engaged 
in such a computation.  Their preliminary work shows that one can take -- unconditionally -- 
$N_0 = 5.5 \times 10^6$.  This $N_0$ is small enough to allow case-by-case analysis, and we believe that the final result will be more like 
$N_0 \approx 500$.  The work will appear elsewhere. 

\subsection{Dramatis Personae and Acknowledgments} \textbf{} \\ \\
The 2007-2008 UGA VIGRE research group in number theory included:
\\ \\
Group leaders (year long): \\ Pete L. Clark (assistant professor), Patrick Corn (postdoc) \\ 
Graduate students (year long): \\ Steve Lane, Jim Stankewicz, Nathan Walters, Steve Winburn, Ben Wyser \\
Graduate students (spring semester only): Brian Cook \\ 
Undergraduate student (year long): Alex Rice.
\\ \\
\noindent
For a 21st century paper on elliptic curves, the theory we need here is relatively middlebrow and classical: most of the 
results we need go back, in some form, to Deuring or even Weber.  Each of the individual results we use can be picked up by a hard-working second year graduate student, but to master them all in a 
limited amount of time while doing research including substantial computer programming is a taller order.  Part of the goal 
of this project was indeed to foster learning by doing, and we have aimed for an exposition which maximizes accessibility to 
the students in the seminar and other early career graduate students.  \\ \indent 
Many of the participants were assigned specific subproblems which they wrote up formally and have been incorporated into this paper.  Specifically, we wish to acknowledge the contributions of Steve Lane in computing 
Table 1, of Alex Rice in $\S 2.4$, of Jim Stankewicz in $\S 5.1$ and of Brian Cook in $\S 8$.  
\\ \indent
The first author would like to thank all the participants in the seminar for an enlightening and stimulating experience; 
this paper represents a substantial advancement of his prior work in this area, which would probably not have been done were 
it not for the interest and involvement of the students.

\section{Background on elliptic curves and complex multiplication}
\noindent

\subsection{Some facts about elliptic curves with complex multiplication}
Let $E$ be an elliptic curve over any field $K$.  A $K$-rational endomorphism of $E$ is a morphism of $K$-varieties 
$\varphi: E \ra E$ such that $\varphi(O) = O$.  Then $\varphi$ induces an endomorphism (i.e., 
self-homomorphism) on the group $E(L)$ of $L$-rational points, for any field extension $L$ of $K$.  
By definition, the \textbf{endomorphism ring} of $E$ is the set of all $\overline{K}$-rational 
endomorphisms of $E$, endowed with the structure of a ring under pointwise addition and composition.  
As for any ring, there is a natural homomorphism $\iota: \Z \ra \End(E)$, in which the image of $n$ is 
the multiplication by $n$ map on $E$, traditionally denoted $[n]$.  
\\ \\
In all cases $\varphi$ is an injection and $\End(E)$, as an abelian group, 
is a free $\Z$-module of rank $1$, $2$ or $4$.  When $\End(E)$ has rank $4$, the endomorphism ring 
is noncommutative, an order in a definite rational quaternion algebra.  Such an elliptic curve is said 
to be \textbf{supersingular}; supersingular elliptic curves over $K$ exist iff $K$ 
has positive characteristic.  So if $K$ has characteristic $0$, we have either $\End(E) = \Z$, or 
$\End(E) \cong \Z^2$ as a free abelian group; in the latter case $\End(E)$ is isomorphic to an order 
$\mathcal{O}$ of an imaginary quadratic field $\Q(\sqrt{-n})$, and ``thus'' we say that 
$E$ has \textbf{complex multiplication}.  More precisely, we say $E$ has $\mathcal{O}$-CM if $\End(E) \cong \mathcal{O}$.  
Since the ring $\mathcal{O}$ has exactly one nontrivial automorphism -- complex conjugation -- if $\End(E) \cong 
\mathcal{O}$, there are two such isomorphisms.  
\\ \\
Let $D_0$ be a \textbf{fundamental imaginary quadratic discriminant}, i.e., the discriminant of the full ring of 
integer of some imaginary quadratic field.  More concretely, $D_0$ is a negative integer which is either (i) congruent 
to $1 \pmod 4$ and squarefree, or (ii) congruent to $0 \pmod 4$ and such that $\frac{D_0}{4}$ is squarefree.  Every 
imaginary quadratic order $\OO$ in $\Q(\sqrt{-D_0})$ is of the form 
$\Z[f\tau_n]$ for a uniquely determined $f \in \Z^+$, the \textbf{conductor} of $\OO$.  Thus an 
order is determined by its fundamental discriminant $D_0$ -- the discriminant of the full ring of integers of 
$\OO \otimes \Q$ -- and $f$.  On the other hand, an order is also determined by its \textbf{discriminant} $D = f^2D_0$.  
This means that for any imaginary quadratic discriminant $D$ -- i.e., an integer $D$ with $D < 0$ and $D \equiv 0,1 \pmod 4$ --
there exists a unique (up to isomorphism) imaginary quadratic order $\OO(D)$ of discriminant $D$.  
\\ \\
For any integral domain $R$, one may consider its \textbf{Picard group} $\Pic(R)$, of rank one locally free $R$-modules 
under tensor product.  Otherwise put, $\Pic(R)$ is the quotient of the group of invertible fractional $R$-ideals by 
the subgroup of principal $R$-ideals.  The \textbf{class number} $h(R)$ is the cardinality of $\Pic(R)$.  For an arbitrary 
domain $R$, the class number may well be infinite, but it is finite when $R$ is an order in any algebraic number field, so 
in particular when $R = R(n,d)$ is an imaginary quadratic order.  When $R$ is a Dedekind domain all nonzero fractional 
ideals are invertible, and $\Pic(R) = \Cl(R)$ is the usual ideal class group.  
\\ \\
We abbreviate $h(\mathcal{O}(D))$ to $h(D)$, and if $K = \Q(D_0)$ is an imaginary quadratic field, then the class 
number of $K$, denoted $h(K)$, means the class number of the maximal order $\OO_K$ of $K$.  
\\ \\
Until further notice we fix an imaginary quadratic order $\mathcal{O}$, of discriminant $D$, and with quotient 
field $K = \Q(\sqrt{D_0})$.    

\begin{fact}
a) There exists at least one complex elliptic curve with $\mathcal{O}$-CM. \\
b) Let $E$, $E'$ be any two complex elliptic curves with $\mathcal{O}$-CM.  The $j$-invariants $j(E)$ and $j(E')$ are 
Galois conjugate algebraic integers.  In other words, $j(E)$ is a root of some monic polynomial with $\Z$-coefficients, and 
if $P(t)$ is the minimal such polynomial, $P(j'(E)) = 0$ also.  \\ 
c) Thus there is a unique irreducible, monic polynomial $H_D(t) \in \Z[t]$ whose roots are the $j$-invariants of the various non-isomorphic $\mathcal{O}$-CM complex elliptic curves. \\
d) The degree of $H_D(t)$ is the class number $h(\mathcal{O}) = h(D)$ of the order $\mathcal{O}$, so when $\mathcal{O}$ 
is the full ring of integers of its quotient field $K$, $\deg(H_D(t)) = h(K)$, the class number of $K$. \\
e) Let $F_D := \Q[t]/H_D(t)$.  Then $F_D$ can be embedded in the real numbers, so in particular is linearly disjoint 
from the imaginary quadratic field $K$.  Let $K_D$ denote the compositum of $F_D$ and $K$.  Then $K_D/K$ is abelian, 
with Galois group canonically isomorphic to $\Pic(\mathcal{O})$.  Moreover, $K_D/\Q$ is Galois and the exact sequence
\[1 \ra \Gal(K_D/K) \ra \Gal(K_D/\Q) \ra \Gal(K/\Q) \ra 1 \]
splits, i.e., $\Gal(K_D/\Q)$ is up to isomorphism the semidirect product of $\Pic(\mathcal{O})$ with the cylic group 
$Z_2$ of order $2$, where the map $Z_2 \ra \Aut(\Pic(\mathcal{O}))$ takes the nontrivial element of $Z_2$ to inversion: 
$x \mapsto x^{-1}$.  
\end{fact}

\noindent
References for this fact include: Cox \cite{Cox} and Silverman II \cite{SilvermanII}.  
\\ \\
This fact has many implications.  First, it follows that one can define an $\mathcal{O}$-CM elliptic curve 
over a number field $F$ iff $F \supset F_D$.  In particular, it follows that one can define an $\mathcal{O}$-CM 
elliptic curve over $\Q$ iff $h(D) = 1$, which by the Heegner-Baker-Stark theorem is known to occur for exactly $13$ 
values of $D$:
\[D = -3, -4, -7, -8, -11, -12, -16, -27, -28, -19, -43, -67, -163. \]
Let $E: y^2 = x^3 + Ax +B$ be a complex elliptic curve in Weierstrass form.  We define a 
\textbf{Weber function} $h$ on $E$, as: \\ \\
$h(x,y) = x$ if $AB \neq 0$, \\ $h(x,y) = x^2$ if $B = 0$, \\
$h(x,y) = x^3$ if $A = 0$. 
\\ \\
(The point of the Weber function is to make explicit the quotient map 
$E \ra E/\Aut(E) \cong \PP^1$.  See \cite[Ch. II]{SilvermanII} for more details.)
\\ \\
If $E$ is defined over some subfield $K$ of $\C$, let $K(E[N])$ be the field extension of $K$ 
obtained by adjoining the coordinates of all the $N$-torsion points on $E$.  
\\ \\
The following is a celebrated classical result.
\begin{thm}(Weber)
Let $D$ be an imaginary quadratic order, and $E_{/F_D}$ an $\mathcal{O}_K$-CM elliptic curve.  For any positive integer 
$N$, the field $\Q(\sqrt{-D},j(E),h(E[N]))$ is the $N$-ray class field of $K = \Q(\sqrt{-D})$.
\end{thm}
\noindent
Proof: See e.g. \cite[Thm. II.5.6]{SilvermanII}.  
%\\ \\
%One reason to take Weber functions is that the extension $K(j(E),E[n])/K$ need not be abelian.  However, $K_D(E[n])/K_D$ 
%is abelian, and this allows us to state a more precise result.  Let $R(N)/K$ be the $N$-ring class field of $K$.  
%
%\begin{thm}
%Let $N$ be an odd prime; $D = f^2 D_0$ be an imaginary quadratic discriminant; and $E_{/K_D}$ an elliptic curve with 
%$\mathcal{O}(D)$-CM.  \\
%a) $K_D(E[N]) = R(Nf)[\zeta_N]$. \\
%b) $[K_D(E[N]):K_D] = \frac{N-1}{w(K)} \cdot \left(N-\left(\frac{D}{N}\right) \right)$. \\
%\end{thm}

\begin{cor}
\label{COR9}
Let $K = \Q(\sqrt{D_0})$ be an imaginary quadratic field, and let $E_{/F(D_0)}$ be an elliptic 
curve with $\mathcal{O}_K$-CM.   Let $N$ be an odd prime.  Then 
\[[\Q(\sqrt{D_0},j(E),h(E[N])):\Q(\sqrt{D_0},j(E))] = 
\left(\frac{N-1}{w(K)}\right)\left(N-\left(\frac{D_0}{N}\right)\right). \]
\end{cor}
\noindent
Proof: We deduce the corollary from the theorem using the description of the $N$-ray class 
field $K(N)$ of $K$ provided by class field theory.  Namely, consider the $N$-\textbf{ring class 
field} $L(N)$, a subextension of $K(N)/K$.  Putting $D = N^2 \cdot D_0$, we have

\[\Gal(L(N)/K) \cong \Pic(\mathcal{O(D)}), \]
whereas 
\[\Gal(K(N)/L(N)) \cong (\Z/N\Z)^{\times}/{\pm 1}. \]
Recall the relative class number formula \cite[Thm. 7.24]{Cox} 
\[ \frac{h(N^2D_0)}{h(D_0)} = \frac{N-\left(\frac{D_0}{N}\right)}{[\mathcal{O}_K^{\times}:\mathcal{O}^{\times}]}, \]
Thus \[ [Q(\sqrt{D_0},j(E),h(E[N]): \Q(\sqrt{D_0},j(E))] = [K(N):K(1)]  \] \[= \frac{[K(N):K]}{[K(1):K]} 
= \frac{ h(N^2D_0) (N-1)}{2h(D_0)} = \frac{N-1}{ w(K)} \cdot \left(N-\left(\frac{D_0}{N}\right)\right). \]

%By arguments similar to those used in the proofs of Theorems \ref{JZEROTHM} and 
%\ref{THEOREM4}, this implies Theorem \ref{NEWTHM} of $\S 1.5$.  

\subsection{The Galois representation}
\noindent
Let $F$ be a field of characteristic $0$, $E_{/F}$ an elliptic curve, and $N$ a positive integer.  Let 
$\sigma \in \Gal_F = \Aut(\overline{F}/F)$.  Let $E[N]$ be the set of $N$-torsion points on $E$ over $\overline{F}$; 
the action of $\Gal_F$ is seen to be $\Z/N\Z$-linear, so $E[N]$ may naturally be viewed as a $\Z/N\Z[\Gal_F]$-module.  
Recall that, as a $\Z/N\Z$-module (or equivalently, as an abelian group), $E[N] \cong \Z/N\Z \times \Z/N\Z$ 
\cite{Silverman}.  It is notationally convenient to choose such an isomorphism -- i.e., to choose an ordered 
$\Z/N\Z$-basis $e_1, \ e_2$ of $E[N]$.  The $\Z/N\Z[\Gal_F]$-module structure is then given by a homomorphism 
\[\rho_N: \Gal_F \ra \GL_2(\Z/N\Z), \]
which we call the \textbf{mod N Galois representation} associated to $E$.  Let $M = F(E[N])$ be the field extension 
obtained by adjoining to $F$ the $x$ and $y$ coordinates of all the $N$-torsion points.  Then the kernel of 
$\rho_N$ is nothing else than $\Gal(\overline{F}/M) = \Gal_M$, so $\rho_N$ factors through to give an embedding 
\[\rho_N: \Gal(M/F) \hookrightarrow \GL_2(\Z/N\Z). \]
There is ``a piece'' of $\rho_N$ which is well understood in all cases.  Namely, composing with the 
determinant map $\det: \GL_2(\Z/N\Z) \ra (\Z/N\Z)^{\times}$, we get a homomorphism
\[ \det(\rho_N): \Gal(M/F) \ra (\Z/N\Z)^{\times}. \]
This homomorphism evidently cuts out an abelian extension of $F$, so can be viewed as a ``character'' of the group 
$\Gal(M/F)$.  More precisely:
\begin{thm}
We have $\det(\rho_N) = \chi_N$, where $\chi_N$ is the mod $N$ \textbf{cyclotomic character}, defined as follows:
\[\chi_N: \Gal_F \ra \Gal(F(\zeta_N)/F) \ra (\Z/N\Z)^{\times}, \]
where $\sigma \in \Gal_F \mapsto \sigma \in \Gal(F(\zeta_N)/F)$, an automorphism which is determined by its effect 
on a primitive $N$th root of unity: \[\zeta_N \mapsto \sigma(\zeta_N) = \zeta_N^{\chi_N(\sigma)}, \]
for a uniquely determined element $\chi_N(\sigma) \in \Z/N\Z^{\times}$.
\end{thm}
\noindent
Proof: See \cite[Ch. III]{Silverman}.
\begin{cor}\label{DETCOR}
We have $\det(\rho_N(\Gal_F)) = 1$ iff $F$ contains the $N$th roots of unity.
\end{cor}
\noindent 
The following is a special case of an extremely important theorem of Serre:
\begin{thm}(Serre's Open Image Theorem, non-CM Case \cite{Serre72})
Let $E$ be an elliptic curve defined over a number field $F$, and suppose that $E$ \textbf{does not} have 
complex multiplication.  \\
a) For all sufficiently large prime numbers $\ell$, $\rho_{\ell}: \Gal_F \ra GL_2(\Z/\ell\Z)$ is surjective. \\
b) There exists a fixed number $B$ such that for all $N \in \Z^+$, 
\[ \# \coker(\rho_N) := \frac{ \# GL_2(\Z/N\Z)}{\# \rho_N(\Gal_F)} \leq B. \]
\end{thm}
\noindent
In other words, part b) says the failure of all the maps $\rho_N$ to be surjective can be measured by a single finite 
quantity.  Since \[GL_2(\Z/\ell_1 \cdots \ell_r \Z) \cong GL_2(\Z/\ell_1 \Z) \times \cdots \times 
GL_2(\Z/\ell_r \Z), \] this in 
fact implies part a).  Note also that we must allow some finite amount of nonsurjectivity, because we are considering an 
elliptic curve $E$ defined over any number field.  So for instance, start with $E$ over $\Q$ and take $F = \Q(E[N])$ 
to be the extension obtained by adjoining all the coordinates of the $N$-torsion points.  For this $E/F$ one tautologically 
has $\rho_N(\Gal_F) = 1$.  Serre himself noted that there is no elliptic curve over $\Q$ for which all the mod $N$ 
Galois representations are surjective.  
\subsection{Galois representation in the CM case} \textbf{} \\ \\
Our interest here is in the fact that this result fails in the presence of CM. 
\\ \\
We assume that $N$ is an \textbf{odd prime}.  
\\ \\
Suppose first that $E/F$ is a $\mathcal{O}(D)$-CM elliptic curve and that $F$ contains the CM field $K = \Q(\sqrt{D})$, so that 
the action of $\mathcal{O}(D)$ is defined and rational over $F$.  Then, in additional to its $\Z/N\Z[\Gal_F]$-module 
structure, $E[N]$ also has the structure of a $\mathcal{O}$-module.  Morever, the $F$-rationality of the endomorphisms 
means precisely that for all $\sigma \in \Gal_F$ and $\varphi \in \mathcal{O}(D)$, we have $\sigma \varphi = \varphi \sigma$, 
i.e., the two actions commute with each other.\footnote{This can be expressed more concisely as the 
fact that $E[N]$ is a $(\Z/N\Z[\Gal_F],\mathcal{O}(D))$-bimodule, but for our purposes there is no particular advantage 
to using this terminology.}  In fact, since $N = 0$ in $E[N]$, $E[N]$ is naturally a $\mathcal{O}(D) \otimes \Z/N\Z = \mathcal{O}(D)/N \mathcal{O}(D)$-module.
\begin{lemma}(\cite[Lemma 1]{Parish})
The $N$-torsion group $E[N]$ is free of rank $1$ as a (right) $\mathcal{O}(D) \otimes \Z/N\Z$-module, i.e., 
isomorphic to $\mathcal{O}(D) \otimes \Z/N\Z$ itself.
\end{lemma}
\noindent
In particular, the natural $\Z/N\Z$-linear action of $\mathcal{O}(D) \otimes \Z/N\Z$ on $E[N]$ is faithful, so we have 
an embedding of $\Z/N\Z$-algebras
\[ \iota: \mathcal{O}(D) \otimes \Z/N\Z \hookrightarrow \End(E[N]) \cong M_2(\Z/N\Z). \]
Let us denote the image of $\iota$ by $C_N$.  Now, for any $\sigma \in \Gal_F$, the matrix $\rho_N(\sigma)$ gives an 
invertible $\mathcal{O}(D) \otimes \Z/N\Z$-linear map of $E[N]$.  Since the $\mathcal{O}(D) \otimes \Z/N\Z$-linear 
endomorphisms of the free one-dimensional module $E[N]$ are precisely multiplication by an element of 
$\mathcal{O}(D) \otimes \Z/N\Z$ and the invertible ones are elements of the unit group of this ring, we conclude 
\[ \rho_N(\Gal_F) \subset C_N^{\times}. \]
This shows that the CM case is much different, because the Galois extension $F(E[N])/F$ is in this case \textbf{abelian} 
and has size at most $\# C_N^{\times}$, or approximately $N^2$, whereas Serre's theorem asserts that in the non-CM 
case $\rho_N(\Gal_F)$ has, for sufficiently large prime $N$, size $\# \GL_2(\Z/N\Z) = (N^2-1)(N^2-N) \sim N^4$.  
\\ \\
To give more precise results, we must consider separately whether $N$ splits, stays inert or ramifies in $\mathcal{O}(D)$.
\\ \\
Case 1 (split case): $(\frac{D}{N}) = 1$.  Then one sees (e.g. by direct computation) that $C_N$, as a $\F_N$-algebra, 
is isomorphic to $\F_N \oplus \F_N$; therefore the unit group $C_N^{\times}$ is isomorphic to 
$(\Z/N\Z)^{\times} \oplus (\Z/N\Z)^{\times}$.  Thus there are precisely two one-dimensional subspaces $V_1$, $V_2$ of 
$E[N]$ which are simultaneous eigenspaces for $C_N$.  By taking generators $e_1$ of $V_1$ and $e_2$ of $V_2$ as 
basis, we get
 
\[ C_N \cong \{ \left[ \begin{array}{cc} a & 0 \\ 0 & b \end{array} \right]  a,b \in \F_N \}.  \]
The same considerations show that there is, up to conjugacy, a unique subalgebra of $M_2(\F_N)$ 
isomorphic to $\F_N \oplus \F_N$; such an algebra is called a \textbf{split Cartan subalgebra} and its unit group a 
\textbf{split Cartan subgroup}.
\\ \\
Case 2 (inert case): $(\frac{D}{N}) = -1$.  Then one sees that $C_N \cong \F_{N^2}$, a finite field of order $N^2$, so 
that $C_N^{\times}$ is cyclic of order $N^2-1$.  Again ones sees that $\F_{N^2}$ is unique up to conjugacy as a subalgebra 
of $M_2(\F_N)$ (e.g. the result is a special case of the Skolem-Noether theorem on simple subalgebras of 
central simple algebras; or just do a direct computation).  Such an algebra is called a \textbf{nonsplit Cartan subalgebra} 
and the unit group is called a \textbf{nonsplit Cartan subgroup}.
\\ \\
Case 3 (ramified case): $N$ divides $D$.  Then $C_N \cong \F_N[t]/(t^2)$, i.e., is generated over the center (the scalar 
matrices) by a single nilpotent matrix $g$.  Since the eigenvalues of $g$ are $\F_N$-rational, we can put $g$ in Jordan 
canonical form, and this gives a choice of basis such that 
\[ C_N \cong \{ \left[ \begin{array}{cc} a & b \\ 0 & a \end{array} \right] \ a,b \in \F_N \} .\]
Again $C_N$ is unique up to conjugacy; for lack of a better name, we shall call it a \textbf{pseudo-Cartan subalgebra}.  
Evidently $C_N \cong Z_{N-1} \oplus Z_N \cong Z_{N^2-N}$.  
\\ \\
We now introduce a third operator on $E[N]$: by Fact 1 above, we can choose an embedding of $K$ into $\C$ which carries $\Q(j_D)$ into the real numbers.  With this understanding, complex conjugation $c$ induces an $\F_N$-linear 
automorphism of $E[N]$.  
\begin{lemma}
Let $N$ be an odd positive integer.  The characteristic polynomial of complex conjugation acting on the free $2$-dimensional $\Z/N\Z$-module $E[N]$ is 
$t^2-1$.
\end{lemma}
\noindent
Proof: Clearly $c$ satisfies the polynomial $t^2-1$, so what we must show is that $c \neq \pm 1$.  If $c = 1$ then 
$c$ acts trivially on each $N$-torsion point and we would have $\dim_{\Z/N\Z} E[N](\R) = 2$.  If $c = -1$ then (since $N$ is odd), 
$c$ acts nontrivially on each $N$-torsion point, and we would have $\dim E[N](\R) = 0$.  But it is easy to see that 
the correct answer is $\dim E[N](\R) = 1$: indeed, a little thought shows that the one-dimensional compact real Lie group 
$E(\R)$ is isomorphic either to $S^1$ (if a defining Weierstrass cubic has one real root) or to $S^1 \times \Z/2\Z$ (if 
all $2$-torsion points (if a defining Weierstrass cubic has three real roots), and either way $E[N](\R) \cong \Z/N\Z$. 
\begin{lemma}(\cite{SerreCM}, \cite{Serre66})
Let $E/\Q(j_D)$ be an $\mathcal{O}(D)$-CM elliptic curve, and let $\sigma$ be the nonidentity element of 
$\Aut(\Q(j_D,\sqrt{D})/\Q(j_D))$. \\ a) As operators on $E[N]$, we have $\sigma = c$. \\
b) Therefore $\Q(j_D,E[N])$ contains $\Q(\sqrt{D})$.
\end{lemma}
\noindent
There is also a natural nontrivial action of complex conjugation on $\mathcal{O}(D)$, and the homomorphism 
$\iota: \mathcal{O}(D) \ra \End(E[N])$ is $c$-equivariant: $\iota \circ c = c \circ \iota$.  This, together with 
the nontriviality of the $c$-action on $\mathcal{O}(D)$, is equivalent to the fact that conjugation by $c$ stabilizes  
$C_N$ and induces a nontrivial involution on it.  
\\ \\
In the split case we find that, with respect to the chosen basis $e_1, \ e_2$ of $C_N$-eigenspaces, $c$ is equal to either
permutation matrix $ \left[ \begin{array}{cc} 0 & 1 \\ 1 & 0 \end{array} \right]$ or its negative.   Either way, 
the effect of conjugation by $c$ is 
$\left[ \begin{array}{cc} a & 0 \\ 0 & b \end{array} \right] \mapsto \left[ \begin{array}{cc} b & 0 \\ 0 & a 
\end{array} \right]$.  Explicit computation shows that the Cartan subgroup $C_N^{\times}$ has index $2$ in its normalizer $N(C_N^{\times})$.
\\ \\
In the inert case, conjugation by $c$ stablizes $C_N \cong \F_{N^2}$ and induces the unique nontrivial Galois 
automorphism, the Frobenius map: $\Frob_N: x \mapsto x^N$.  The elements of $N(C_N^{\times}) \setminus C_N^{\times}$ correspond to $\Frob_N$-semilinear automorphisms 
of the $1$-dimensional $\F_{N^2}$-vector space $V = E[N]$, i.e., maps $\sigma: V \ra V$ such that for $v,w \in V$, 
$\sigma(vw) = \Frob_N(v) \sigma(w)$.  Such a map is uniquely specified by $\sigma(1)$, so that $\# N(C_N^{\times}) \setminus 
C_N^{\times} = N^2-1$, i.e., $[N(C_N^{\times}):C_N^{\times}] = 2$.
\\ \\
In the ramified case, complex conjugation induces a nontrivial involution of the (non-semisimple) $\F_N$-algebra 
$C_N \cong \F_N[t]/(t^2)$.  The automorphism group $\Aut(C_N/\F_N)$ is isomorphic to $\cong F_{N-1}^{\times}$ so 
has a unique element of order $2$, $t \mapsto -t$.  Therefore conjugation by $c$ has the effect 
$ \left[ \begin{array}{cc} a & b \\ 0 & a \end{array} \right] \mapsto \left[ \begin{array}{cc} a & -b \\ 0 & a 
\end{array} \right].$  Note that this case is different from the previous two in that the normalizer of $C_N^{\times}$ 
is the entire Borel subgroup $\{\left[ \begin{array}{cc} a & b \\ 0 & c \end{array} \right] \ | \ a,b,c \in \F_N, ac \neq 0 \}$.
\\ \\
Given all this information, one readily deduces the following result:
\begin{thm}
Let $F$ be a number field, and $E_{/F}$ an elliptic curve with $\mathcal{O}(D)$-CM.  Let $M = F(E[N])$ be the field 
extension of $F$ obtained by adjoining $x$ and $y$ coordinates of all the $N$-torsion points of $E$.  \\ 
a) The CM field $K = \Q(\sqrt{-D})$ is contained in $M$, so we get a short exact sequence 
\begin{equation}
\label{EQONE}
1 \ra \Gal(M/KF) \ra \Gal(M/F) \ra \Gal(KF/F) \ra 1. 
\end{equation}
b) Under the natural embedding $\rho_N: \Gal(M/F) \hookrightarrow GL_2(\F_N)$, the subgroup $\Gal(M/KF)$ 
embeds in the unit group $C_N^{\times}$.  \\
c) The sequence (\ref{EQONE}) splits, with a splitting given by a choice of an involution $c \in 
N(C^{\times}) \setminus C^{\times}$.
\end{thm}
\noindent
This result gives upper bounds on the the degree $[F(E[N]):F]$ which improve upon the obvious bound of $\# \GL_2(\F_N)$: 
\begin{cor}
a) If $(\frac{D}{N}) = 1$, then $[F(E[N]):F] \ | \ 2 (N-1)^2$. \\
b) If $(\frac{D}{N}) = -1$, then $[F(E[N]):F] \ | \ 2 (N^2-1)$. \\
c) If $(\frac{D}{N}) = 0$, then $[F(E[N]):F] \ | \ 2(N^2-N)$.
\end{cor}
\noindent
Proof: Using the exact sequence (\ref{EQONE}) we see that 
\[\# \Gal(M/F) = \# \Gal(M/KF) \cdot \# \Gal(K/F) \ | \ \# (C_N)^{\times} \cdot 2. \]
And we know that $C_N^{\times}$ has order $(N-1)^2$, $(N^2-1)$ or $N^2-N$ according 
to whether $N$ splits, is inert, or is ramified in $\mathcal{O}(D)$.  
\\ \\
The slogan here is that the image of the Galois representation $\rho_N$ should be ``as large as possible'', up to a factor which 
is uniformly bounded as $N$ varies, but in the CM case $GL_2(\F_N)$ is impossibly large.  The correct answer is again due to 
Serre:
\begin{thm}(Open Image Theorem, CM case \cite{Serre66}):
\label{SERRE1}
Let $F$ be a number field and $E_{/F}$ be an elliptic curve with $\mathcal{O}$-CM.  Then for all sufficiently 
large primes $N$, we have: \\
$\bullet$ $\rho_N(\Gal_F)) = N(C_N)$, if $K = \Q(\sqrt{-D})$ is not contained in $F$, \\
$\bullet$ $\rho_N(\Gal_F) = C_N^{\times}$, if $K \subset F$.
\end{thm}
\noindent
Since Serre's theorem only holds for sufficiently large primes $N$, the case of $N \ | \ D$ can be 
completely ignored.  Nevertheless Theorem \ref{SERRE1} tells us to ``expect'' that the $N$-torsion fields will be 
as large as possible.  In the next section we use elementary group theory to deduce consequence for the least 
degree of an $N$-torsion point.
\subsection{Orbits under $C_N^{\times}$ and applications} \textbf{} \\ \\
We maintain the notation of the previous section: $E_{/F}$ is an elliptic curve with $\mathcal{O}(D)$-CM; 
$N$ is an odd prime number; $C_N = \iota(\mathcal{O} \otimes \Z/N\Z) \subset \End(E[N])$; $C_N^{\times}$ 
is the unit group of $C_N$; $N(C_N^{\times})$ is the normalizer.  
\begin{lemma}
\label{ORBITLEMMA}
a) The orbits of $C_N^{\times}$ on $E[N] \setminus \{0\}$ are as follows: \\
(i) If $(\frac{D}{N}) = 1$, the two one-dimensional eigenspaces for $C_N$ give two orbits of size 
$N-1$; all the remaining points lie in a single orbit of size $(N-1)^2$.  \\
(ii) If $(\frac{D}{N}) = -1$, $E[N] \setminus \{0\}$ forms a single $C_N^{\times}$-orbit. \\
(iii) If $(\frac{D}{N}) = 0$, the unique one-dimensional eigenspace for $C_N$ gives an orbit of 
size $N-1$; the remaining points form a single orbit of size $N^2-N$. \\
b) If $(\frac{D}{N}) = 1$, the two orbits of size $N-1$ for $C_N^{\times}$ form a single orbit for 
$N(C_N^{\times})$.
\end{lemma}
\noindent
Proof: A pleasant elementary computation that we leave to the reader.
\\ \\
In the statement of the following result we employ the following convention: if $p$ and $q$ 
are nonzero rational numbers, we say $p \ | \ q$ if $\frac{q}{p} \in \Z$.
\begin{cor}
\label{ORBITCOR}
Let $E_{/F}$ be an $\mathcal{O}(D)$-CM elliptic curve defined over a number field $F$.  Suppose that the image 
$\rho_N(\Gal_{KF})$ of the mod $N$ Galois representation has index $I$ in 
$C_N^{\times}$.  Let $P \in E(\C)$ be any point of exact order $N$, and let $F(P)$ be the extension of $F$ obtained 
by adjoining the coordinates of $P$. \\
(i) If $(\frac{D}{N}) = 1$ and $\sqrt{D} \in F$, then $\frac{1}{I} (N-1) \ | \ [F(P):F]$ \\
(ii) If $(\frac{D}{N}) = 1$ and $\sqrt{D}$ is not in $F$, then $\frac{2}{I} \ | \ [F(P):F]$. \\
(iii) If $(\frac{D}{N}) = -1$, then $\frac{1}{I}(N^2-1) \ | \ [F(P):F]$. \\
(iv) If $(\frac{D}{N}) = 0$, then $\frac{1}{I} (N-1) \ | \ [F(P):F]$.
\end{cor}
\noindent
Proof: Consider of field extensions $F \subset F(P) \subset F(E[N])$.  Then $F(E[N])/F(P)$ is Galois, 
with Galois group canonically isomorphic to $\rho_N(\Gal_F) \cap G(P)$, where $G(P) \subset \GL_2(\F_N)$ is the stabilizer 
of the point $P$.  By the orbit-stabilizer theorem, $[F(P):F]$ is equal to the orbit of $P$ under the action of 
$\Gal_F$.  \\ \indent In case (i) we have $\sqrt{D} \in F$, so that the image of Galois lies in the split Cartan subgroup 
$C_N^{\times} \cong \F_N^{\times} \oplus \F_N^{\times}$.  By Lemma \ref{ORBITLEMMA} the full $C_N^{\times}$-orbits have 
sizes $N-1$ and $(N-1)^2$.  Since we are assuming that $[C_N^{\times}:\rho_N(\Gal_F)] \ | \ I$, it follows that every 
$\rho_N(\Gal_F)$-orbit has size a multiple of $\frac{N-1}{I}$.  Case (ii) is similar except in this case replace the 
gcd of all sizes of $C_N^{\times}$ orbits with the gcd of all sizes of $N(C_N^{\times})$-orbits, which according to 
Lemma \ref{ORBITLEMMA} is $2(N-1)$.  Parts (iii) and (iv) are similar, except here it does not matter whether 
$\sqrt{D}$ lies in the ground field $F$: in case (iii) this is because the orbit size for $C_N^{\times}$ is already 
as large as possible; in case (iv) this is because the minimal $C_N^{\times}$-orbit is stable under complex conjugation.
\section{Proof of Theorem \ref{JZEROTHM}}
\noindent
As in Remark 1.2, Theorem \ref{JZEROTHM}a) is precisely the $D = -3$ case of Theorem \ref{THEOREM4}a).  Indeed, for 
$D = -3$, $w(D) = 6$, and an odd prime splits completely in $\Q(\sqrt{-3})$ iff $N \equiv 1 \pmod 3$.  
\\ \indent
Now suppose we have an $\OO(D)$-CM point on $X_1(N)$ of degree $D$.  If $D = -3$, then according to 
Theorem \ref{THEOREM4}, if $N$ is greater than or equal to some absolute constant $N_1$, we have $d \geq \frac{N-1}{3}$ if 
$N \equiv 1 \pmod 3$ and $d \geq \frac{N^2-1}{6}$ if $N \equiv -1 \pmod 3$. 
\\ \\
The second case is $D = -4$, so $w(D) = 4$, and then Theorem \ref{THEOREM4} says that for $N$ greater than or equal to 
another absolute constant $N_2$, we have $d \geq \frac{N-1}{2}$ if $N \equiv 1 \pmod 4$ and $d \geq \frac{N^2-1}{4}$ if 
$N \equiv -1 \pmod 4$.  
\\ \\
The third case is any other $D$, so $w(D) = 2$ and then by Theorem \ref{SPYTHM}, $d \geq \frac{N-1}{2}$.  Altogether we see that if $N \geq \max(5,N_1,N_2)$ then $d \geq \frac{N-1}{3}$ in all cases, equality 
can be met iff $N \equiv 1 \pmod 3$ (necessarily for an $\OO(-3)$-CM elliptic curve of $j$-invariant $0$), and the next 
smallest possible degree is $\frac{N-1}{2}$, for an $\OO(-4)$-CM elliptic curve of $j$-invariant $1728$.  This completes 
the proof of Theorem \ref{JZEROTHM}.

\section{Proof of Theorem \ref{NEWTHM}} \noindent
Let $N$ be an odd prime number; let $K = \Q(\sqrt{D_0})$ be an imaginary quadratic field; and let 
$E_{/F}$ be an $\mathcal{O}_K$-CM elliptic curve.  Suppose that there exists a point $P \in E(F)$ of 
order $N$.  Let $M$ be the compositum of the CM field $K$ with the $N$-torsion field $F(E[N])$.  We know 
that $K(j(E)) = K(1)$ is the Hilbert class field of $K$ and $M/K(N)$ is abelian of degree divisible by 
$\frac{N-1}{w(K)} \cdot (N-\left(\frac{D_0}{N} \right))$ by Corollary \ref{COR9}.
\\ \\
Split case: $(\frac{D_0}{N}) = 1$.  We know that $\Gal(M/K(1))$ is contained in a split Cartan 
subgroup $C(N) \cong (\Z/N-1\Z)^2$ of $\GL_2(\Z/N\Z)$ with index dividing $w(K)$.  If we had equality -- i.e., 
$[M:K(1)] = (N-1)^2$ --
then by the work of the previous section, for any $N$-torsion point $P \in E(\C)$ we must have 
$N-1 \ | \ [K(1)(P):K]$.  Moreover, as we saw above, passing to a subgroup of index $i$ cuts down this degree by at 
most a factor of $i$, so 
\[ \frac{N-1}{w(K)} \ | \ [K(1)(P):K(1)], \]
and therefore
\[ \frac{h(K)}{w(K)} \cdot (N-1) \ | \ [K(1)(P):K] \ | \ [KF:K]. \]
Since 
\[ \delta(F,K) \cdot [KF:K] = [F:\Q], \]
Theorem \ref{NEWTHM}a) follows. 
\\ \\
Ramified case: $(\frac{D_0}{N}) = 0$.  In this case $\Gal(M/K(1))$ is contained in a pseudo-Cartan 
subgroup $C_N^{\times} \cong \Z/N\Z \times \Z/N-1$ of $\GL_2(\Z/N\Z)$ with index dividing $w(K)$.  As above, the smallest 
orbit of $C_N^{\times}$ on the $N$-torsion has size $N-1$, leading to the bound 
\[ \frac{h(K)}{w(K)} \cdot (N-1) \ | \ [F:\Q]. \]
In this case, $\Gal(KF/K)$ acts trivially on the unique $N$-torsion subgroup stabilized by $C_N^{\times}$.  From this, one sees that 
we gain an extra factor of $2$ iff $K$ \emph{does} contain $F$, giving the divisibility relation as in Theorem 
\ref{NEWTHM}b).
\\ \\
Inert case: $(\frac{D_0}{N}) = -1$.  In this case $\Gal(M/K(1))$ is contained in a nonsplit Cartan 
subgroup $C_N^{\times} \cong (\Z/(N^2-1)\Z)$ of $\GL_2(\Z/N\Z)$ with index dividing $w(K)$.  As above, 
the $N$-torsion points form a single orbit under $C_N^{\times}$, so arguing as in the split case we get 
\[ \frac{h(K)}{w(K)} \cdot (N^2-1) \ | \ [F:\Q]. \]
This completes the proof of Theorem \ref{NEWTHM}.

\section{Proof of Theorem \ref{THEOREM4}}

\subsection{A technical lemma} \textbf{} \\ \\
Let $w$ be a positive even integer, and let $\zeta = \zeta_{w} = e^{2\pi i /w}$ be a primitive $w$th root of unity.  
Let $G = \langle s \ | \ s^w = 1 \rangle$ be a cyclic group of 
order $w$.  Let $M$ be an abelian group endowed with the following additional structures: \\ \\
$\bullet$ a $\Z$-linear action of $G$, and \\
$\bullet$ A ring homomorphism $\Z[\zeta] \ra \End(M)$. \\ \\
We require first that $\zeta^{\frac{w}{2}} \cdot x = -x$ for all $x \in M$.  We also require that these two actions commute 
with each other: for all $x \in M$, $\zeta \sigma x = \sigma \zeta x$.  
\\ \\
For $i \in \Z/w\Z$, we define $M_i = \{x \in M \ | \ \sigma x = \zeta^i x \}$, and 
\[\textbf{M} = \bigoplus_{i \in \ \Z/w\Z} M_i. \]
Consider the $\Z$-module homomorphism $\Phi: \textbf{M} \ra M$ given 
$(x_i) \mapsto \sum_i x_i$.  Let $\tilde{\Phi} = \Phi \otimes_{\Z} \Z[\frac{1}{w}]: 
\textbf{M}' = \textbf{M} \otimes \Z[\frac{1}{w}] \ra M' = M \otimes \Z[\frac{1}{w}]$.
\begin{lemma}
\label{TECHLEMMAJIM}
Both $\ker(\Phi)$ and $\coker(\Phi)$ are $w$-torsion $\Z$-modules.  It follows that: \\
a) The map $\tilde{\Phi}$ is an isomorphism of $\Z[\frac{1}{w}]$-modules. \\
b) We have $\dim_{\Q} (\textbf{M} \otimes \Q) = \dim_{\Q} (M \otimes \Q)$, and for any 
prime $p$ not dividing $w$, $\Phi$ induces an isomorphism from the 
$p$-primary torsion subgroup $\textbf{M}[p^{\infty}]$ of $\textbf{M}$ to the $p$-primary 
torsion subgroup $M[p^{\infty}]$ of $M$.
\end{lemma}
\noindent
Proof: It is enough to show that the kernel and cokernel of $\Phi$ are $w$-torsion; for if so, tensoring the short 
exact sequences
\[0 \ra \Ker(\Phi) \ra \textbf{M} \stackrel{\Phi}{\ra} \Phi(\textbf{M}) \ra 0\]
and 
\[0 \ra \textbf{M}/\Ker(\Phi) \stackrel{\Phi}{\ra} M \ra \Coker(\Phi) \ra 0 \]
of $\Z$-modules with the flat $\Z$-module $\Z[\frac{1}{w}]$ shows that $\tilde{\Phi}$ is an isomorphism.
\\ \\
Step 1: We show $\Ker(\Phi) = \Ker(\Phi)[w]$.  Let
$P =(P_0,\dots,P_{w-1})$ be an element of $\ker{\Phi}$, so that $$P_0 + \dots + P_{w-1} = 0.$$  Applying $\sigma$, 
we obtain
$$P_0 + \zeta P_1 + \dots + \zeta^{w-1}P_{w-1} = 0.$$
Applying $\sigma$ $w-2$ more times, we arrive at the matrix equation $AP = 0$, where $$A = \left(\begin{array}{ccccc} 1 & 1 & \dots & 1 \\ 1 & \zeta & \dots & \zeta^{w-1} \\ \vdots & \vdots & \ddots & \vdots \\ 1 & \zeta^{w-1} & \dots & \zeta^{(w-1)(w-1)}\end{array}\right).$$

It is therefore also a solution to $A^2 P = 0$, where $$A^2 = \left(\begin{array}{ccccc} w & 0 & \dots& \dots & 0 \\ 0 & 0 & \dots & 0 & w \\ \vdots & \vdots & \adots & w & 0 \\ \vdots & 0 & \adots & \adots & \vdots \\ 0 & w & 0 & \dots & 0 \end{array}\right).$$
Thus $wP_0 = wP_{w-1} = \dots = wP_1 = 0$, i.e., $w P = 0$.  
\\ \\
Step 2: We show $\Coker(\Phi) = \Coker(\Phi)[w]$.  Let $P \in M$.  Define a $w \times w$ matrix 
$$B =\left(\begin{array}{ccccc} P  & \sigma(P)  & \sigma^2(P)  & \dots  & \sigma^{w-1}(P)\\ P  & \zeta^{-1}\sigma(P)  & \zeta^{-2}\sigma^2(P)  & \dots  &\zeta^{-(w-1)}\sigma^{w-1}(P)\\ \vdots  & \vdots  & \vdots & \dots & \vdots \\ P  & \zeta^{-(w-1)}\sigma(P)  & \zeta^{-2(w-1)}\sigma^2(P) & \dots  & \zeta^{-(w-1)(w-1)}\sigma^{w-1}(P) \end{array}\right).$$
Notice that the sum of all the entries of $B$ is $wP$: indeed, this is the sum of the entries in the first column, and 
since for any $j \neq 0 \pmod w$ we have $\sum_{i=1}^{w-1} \zeta^{-ji} = 0$, each of the other columns sums to $0$.  Now 
for $1 \leq i \leq w$, put \[P_{i-1} = \sum_{k = 0}^{w-1} \zeta^{-ki}\sigma^k(P). \]  Then 
\[ \sigma(P_i) = \sum_{k = 0}^{w-1} \zeta^{-ki}\sigma^{k+1}(P) = 
\zeta^i\sum_{k = 0}^{w-1} \zeta^{-(k+1)i}\sigma^{k+1}(P) = \zeta^i P, \]
so $P_{i-1} \in M_i$.  Therefore 
\[wP = \Phi((P_0,\ldots,P_{w-1})) \in \Phi(\textbf{M}). \]
This completes the proof of the lemma.

\subsection{Application to the proof of Theorem \ref{THEOREM4}} \textbf{} \\ \\
Now let $\mathcal{O}$ be an imaginary quadratic order of discriminant $D$, $K = \Q(\sqrt{D})$, and let
$N > w = w(\mathcal{O})$ be a prime which splits in $K$.  Let $K_D = K(j_D)$, and let $E/K_D$ be an 
$\mathcal{O}$-CM elliptic curve.  By the work of $\S 2.3$, we know that there 
exists an extension $K_D(P)/K_D$, which is cyclic of degree dividing $N-1$, such that over $K_D(P)$ 
$E$ has a point $P$ of exact order $N$.  Let us first assume that $[K_D(P):K_D] = N-1$; afterwards we 
will discuss how to modify the argument to deal with the case in which the degree strictly 
divides $N-1$. \\ \indent
Our assumptions imply that $N \equiv 1 \pmod w$.  Therefore, by Galois theory, there exists a unique 
subextension $K_D \subset L \subset K_D(P)$ with $G = \Gal(K_D(P)/L)$ cyclic of order $w$.  Now we are in 
the setup of the previous section: take $M = E(K_D(P))$; the $G$-action is the restriction of 
the natural $\Gal(K_D/K_D)$-action on $E(K_D(P))$, the $\Z[\zeta]$-action comes from the fact that 
$\mathcal{O} = \End(E)$ contains the $w$th roots of unity, and the compatibility of these two actions 
is a consequence of the rationality of the endomorphisms over $K_D$ (hence also over 
$L$).  Since $E(K_D(P))$ contains a point whose order is a prime $N$ not divisible by $w$, by Lemma 
\ref{TECHLEMMAJIM} there exists some $i \in \Z/w\Z$ such that $M_i$ contains an element of order $N$.  \\ \indent
Using the theory of twisting in the Galois cohomology of elliptic curves, we may interpret 
$M_i$ as the group of $L$-rational points on a $K_D(P)/L$-twisted form of the elliptic curve $E$.  
Specifically, the set of such twisted forms are parameterized by 
\[H^1(\Gal(K_D(P)/L),\Aut(E)) = \Hom(G,\Z/w\Z) \cong \Z/w\Z, \]
the last isomorphism being given by \[(\varphi: G \ra \Z/w\Z) \mapsto \zeta^{i_{\varphi}} = 
\varphi(\sigma). \]
Corresponding to $\zeta^i = \zeta^{i_{\varphi}} \in \Z/w\Z$ we build a twisted $\Gal(K_D(P)/L)$-action on 
$E(K_D(P))$: 
\[\sigma \cdot_i x := \zeta^{-i} \sigma x. \]
This is exactly the relation defining $M_i$.  In other words, the abstract decomposition of 
the $\Z[\frac{1}{w}]$-module $\textbf{M}' \stackrel{\sim}{\ra} M'$ corresponds to a decomposition 
of the Mordell-Weil group -- up to $w$-torsion -- of $E(K_D(P))$ into a direct sum of the 
Mordell-Weil groups of the $w$ different twists of $E_{/L}$ via the cyclic extension $K_D(P)/L$ and 
the automorphism group of $E$.  (When $w = 2$, this result -- decomposition of the Mordell-Weil 
group under a quadratic extension -- is very well known.)  Thus we have produced an $\OO$-CM 
elliptic curve over a field of degree $\frac{2(N-1)}{w(\OO)}$ with a rational $N$-torsion point, giving 
the statement of Theorem \ref{THEOREM4}a). \\ \indent
It remains to deal with the case in which $d = [K_D(P):K]$ strictly divides $N-1$.  If $w \ | \ d$, we 
can run through the above argument verbatim, getting in fact an $\OO$-CM elliptic curve with a rational 
$N$-torsion point over a field of degree $\frac{2d}{w}$, which is \emph{a priori} stronger than what 
we are trying to prove.  This necessarily is the case if $w = 2$.  If $w = 4$ and $d$ is a multiple of 
$2$ but not a multiple of $4$, we run through the above argument using quadratic twists instead of 
quartic twists.  If $w = 6$ and $d$ is a multiple of $2$ but not of $6$, then we run through the 
above using quadratic twists instead of sextic twists.  One sees easily that we get exactly 
the same bounds.  This completes the proof of Theorem \ref{THEOREM4}a).  
\\ \\
Proof of Theorem \ref{THEOREM4}b): Suppose first that $N$ is an odd prime with $(\frac{D}{N}) = 1$.  Let 
$F_D = \Q(j(E)) = \Q(j_D)$ be the number field generated by the j-invariant of the quadratic order 
$\mathcal{O}(D)$, and let $E_{/F_D}$ be any $\mathcal{O}(D)$-CM elliptic curve.  
Serre's Theorem \ref{SERRE1} says that there exists $N_0 = N_0(D)$ such that if $N \geq N_0$, the image 
$\rho_N(\Gal_{F_D})$ in $\GL_2(\F_N)$ will be $N(C_N^{\times})$, the normalizer of a split Cartan subgroup, 
and then Corollary \ref{ORBITCOR} applies to show that the least degree $[F_D(P):F_D]$ is a multiple of $2(N-1)$.   
\\ \indent
Now suppose that we have any number field $F$, $E'_{/F}$ an $\mathcal{O}(D)$-CM elliptic curve with an $F$-rational 
point of prime order $N \geq N_0$.  The theory of twisting -- together with the Kummer isomorphism $H^1(\Gal_F,\mu_d) 
\cong F^{\times}/F^{\times d}$ -- implies first that $F \supset F_D$, and second that there exists an extension $L$ of 
$F$, of degree $w(\mathcal{O})$ such that $E_{/L} \cong E'_{/L}$.  Therefore, since $E'$ has an $F$-rational torsion point 
of order $N$, $E$ has an $L$-rational torsion point of order $N$, so 
\[ 2(N-1) \ | \ [F_D:\Q] \ | \ [L:F_D][F_D:\Q] = [L:\Q] = [L:F][F:\Q] = w(\mathcal{O})[F:\Q], \]
and hence \[\frac{2(N-1)}{w(\mathcal{O})} \ | \ [F:\Q]. \]
The argument in the case $(\frac{D}{N}) = -1$ is quite similar: then there exists $N_0$ such that $N \geq N_0$ implies 
that, for our fixed $E_{/F_D}$ as above  we have $[F_D(P):F_D] = N^2-1$ (note that this is the order of the stabilizer of $P$ 
in all of $GL_2(\F_N)$, hence the largest possible order, so there is no further contribution coming from the action of 
complex conjugation) and arguing as before we get \[\frac{N^2-1}{w(\mathcal{O})} \ | \ [F:\Q]. \]
Since we are taking $N$ arbitrarily large compared to $D$, we do not have to worry about the ramified case. 

\section{Proof of Theorem \ref{LEASTDEGTHM}} 
\noindent
For a negative quadratic discriminant $D$, write $d_D(N)$ for the least degree of an $\mathcal{O}(D)$-CM point on 
$X_1(N)$, and $d_{\CM}(N)$ for the least degree of a CM point on $X_1(N)$, so $d_{\CM}(N) = \min_{D} d_D(N)$.  
\\ \\
We will need the following two estimates:
\begin{lemma}
\label{LASTLEMMA}
Suppose $D$ is a positive integer and $N$ a prime, with $(\frac{-D}{N}) = 1$.  Then there exists 
a CM point on $X_1(N)$ of degree dividing $2(N-1)h(\Q(\sqrt{-D}))$.
\end{lemma}
\noindent 
Proof of Lemma \ref{LASTLEMMA}: this is an immediate consequence of the theory of Galois representations 
on CM elliptic curves as recalled in $\S 2.3$.
\begin{lemma}\label{CLASSICLEMMA}
As $D$ tends to $-\infty$ through quadratic discriminants (i.e., $D \equiv 0,1 \pmod 4$), 
the class number $h(D)$ of the imaginary order of discriminant $D$ is $O(\sqrt{D} \log D)$.
\end{lemma}
\noindent
Proof: A consequence of Dirichlet's class number formula; see e.g. \cite[$\S$ 4.2]{Cohen}.
\\ \\
Proof of Theorem \ref{LEASTDEGTHM}: If $N \equiv 1 \pmod 4$, by Theorem \ref{THEOREM4} we 
have $d_{\CM}(N) \leq \frac{N-1}{2}$.  This is stronger than the bounds we are claiming for 
arbitrary $N$, so we may assume that $N \equiv -1 \pmod 4$.  
\\ \indent
For such $N$, let $D$ be a negative quadratic discriminant not divisible by $N$.  Then
\[ 1 = \left(\frac{D}{N}\right) \iff \left(\frac{|D|}{N}\right) = -1, \]
so we are interested in the least positive integer $M$ which is first, a quadratic nonresidue modulo $N$ 
and second, is congruent to $0$ or $-1$ modulo $4$, so that $-M$ is an imaginary quadratic discriminant. \\ 
\indent  
In fact this latter condition is nothing to worry about: let $M$ be the least positive quadratic 
nonresidue modulo $N$.  Then certainly $M$ is squarefree, so $M$ is not $0 \pmod 4$.  
If $M \equiv -1 \pmod 4$, then $D = -M$ is the discriminant of $\Q(\sqrt{-M})$.  
If $M \equiv 1,2 \pmod 4$, then it is not $-M$ but $-4M$ which is the discriminant of $\Q(\sqrt{-M})$.  
But if $M$ is a quadratic nonresidue modulo the odd prime $N$, so is $4M$, and if we know that 
$M = O(f(N))$ for some function $f$, then of course the same holds for $4M$.
\\ \indent
So what is the order of the least quadratic nonresidue modulo $N$?  This is a 
famous classical problem.  The trivial bound -- taking into account only that there are in all 
$\frac{N-1}{2}$ quadratic nonresidues -- is $\frac{N}{2}$, but a bit of thought and experimentation 
suggests that $M$ should be considerably smaller than this.  Long ago Vinogradov conjectured that 
$M = O_{\epsilon}(N^{\epsilon})$, i.e., that $M$ grows more slowly than any power of $N$, but we are still 
far away from an unconditional proof of this.  In 1952 N.C. Ankeny showed that, \textbf{conditionally} on 
GRH, $M = O((\log N)^2)$ \cite{Ankeny}.  In his review of this paper \cite{Erdos}, P. Erd\"os remarks that it is known 
that $M$ is \emph{not} $O(\log N)$, so that Ankeny's bound seems to get admirably close to the truth.  
Vinogradov himself was able to show unconditionally that $M = o(N)$; for more than fifty years, the 
best unconditional bound has been due to D.A. Burgess: $M = O_{\epsilon}(N^{c+\epsilon})$, 
where $c = \frac{e^{-1/2}}{4} = 0.15\ldots$ is ``Burgess' constant'' \cite{Burgess}.  
\\ \\
So, for a large prime $N$, let $M$ be the least quadratic nonresidue modulo $N$ and 
$D = -M$ if $M \equiv -1 \pmod 4$ and $D = -4M$ otherwise.  Applying Lemma \ref{LASTLEMMA} and then 
Lemma \ref{CLASSICLEMMA}, we get
\[d_{\CM}(N) =  O(N h(D)) = O(N \sqrt{|D|} \log |D|). \]
Substituting in the unconditional Burgess bound for $D$, we get 
\[d_{\CM}(N) = O_{\epsilon}(N^{1+c/2 + \epsilon/2} \log(N^{c+\epsilon})). \]
That this bound hold for all $\epsilon > 0$ is equivalent to 
\[d_{\CM}(N) = O_{\epsilon}(N^{1+c/2 + \epsilon}). \]
Applying instead Ankeny's bound, we get, conditionally on GRH,
\[d_{\CM}(N) = O(N \sqrt{(\log N)^2} \log (\log N)^2 = O(N \log N \log \log N). \]

\section{Proof of Theorem \ref{NONLINEARTHM}} \noindent
Although not necessary from a logical point of view, we believe it will make 
for easier reading if we discuss first the special case in which the endomorphism ring is the maximal order and 
second the (less) special case in which the conductor of the order is prime to $N$ before discussing the general case.
\\ \\
\textbf{Case 1: fundamental discriminants} Suppose that there exists some positive number $C$ such that for every odd prime $N$, there exists a 
point on $X_1(N)$ with CM by the \emph{full} ring of integers of some imaginary quadratic field, and of degree at most 
$CN$.  We will derive a contradiction.
\\ \indent
$H := 6C+1$.  Recall that the set of negative quadratic discriminants $D$ such that 
$h(D) \leq H$ is finite \cite{Deuring}, \cite{Heilbronn}, \cite{Siegel}.
Let us write out this set 
as $\{D_1,\ldots,D_n \}$.  
\\ \indent
Let $\mathcal{P}_1$  be the set of primes which are $1 \pmod 4$ and divide $D_k$ for some $1 \leq k \leq n$.  Put 
$R = \# \mathcal{P}_1$.  Similarly, let $\mathcal{P}_3$ be the set of primes which are $3 \pmod 4$ and divide some $D_k$.  Put 
$S = \# \mathcal{P}_3$.
\begin{lemma}
The set $\mathcal{P}_H$ of odd primes $N$ such that 
$\{ (\frac{D}{N}) = -1 \  \forall D \ | \ h(D) \leq H \}$ is infinite; indeed it has density at 
least $(\frac{1}{2})^{R+S+2}$. 
\end{lemma}
\noindent
Proof: Let $N$ be any prime number satisfying: \\
(i) $N \equiv 7 \pmod 8$; \\
(ii) $(\frac{N}{p}) = 1$ for all $p \in \mathcal{P}_1$. \\
(iii) $(\frac{N}{q}) = -1$ for all $q \in \mathcal{P}_3$.  \\ \\
By the Cebotarev density theorem (or even the quantitative version of Dirichlet's theorem on primes in arithmetic 
progressions), the set of such primes $N$ has density $(\frac{1}{2})^{R+S+2}$.  We claim that all such primes lie 
in $\mathcal{P}_H$.  Indeed, we may write 
\[D_k = (- 1) \cdot 2^{a+2b} p_1 \cdots p_r q_1 \cdots q_s = (-1)^{s+1} 2^{a+2b} \prod_{i=1}^r p_i \prod_{j=1}^s (-q_j), \]
where $a,b \in \{0,1\}$, $p_i \in \mathcal{P}_1$ and $q_j \in \mathcal{P}_3$.  Then 
\[\left(\frac{D_k}{N} \right) = \left(\frac{-1}{N} \right)^{s+1} \left(\frac{2}{N} \right)^{a+2b} 
\prod_{i=1}^r \left(\frac{p_i}{N}\right) 
\prod_{j=1}^s \left(\frac{-q_j}{N}\right) = \]
\[ (-1)^{s+1} \cdot 1 \cdot \prod_{i=1}^r \left(\frac{N}{p_i}\right) \cdot \prod_{j=1}^s 
\left(\frac{N}{q_j}\right) = (-1)^{s+1} (-1)^s = -1. \] 
Let $N > H$ be a prime in $\mathcal{P}_H$, and let $D$ be any negative quadratic discriminant.  If $(\frac{D}{N}) = -1$, 
then by Theorem \ref{NEWTHM} we have $d_{D}(N) \geq \frac{N^2-1}{6}$, which for sufficiently large $N$, is 
greater than $CN$.  Otherwise $(\frac{D}{N}) \neq -1$, and by Theorem \ref{NEWTHM} we have 
\[d_{D}(N) \geq \frac{h(D)}{6} (N-1) > \frac{H}{6} (N-1) > CN, \]
since $N > H$.  
\\ \\
\textbf{Case 2: Orders of conductor prime to $N$}: Suppose that $\mathcal{O}(D)$ is an order of conductor $f$ in the 
imaginary quadratic field $K = \Q(\sqrt{D_0})$; let $F$ be a number field and $E_{/F}$ be a 
$\mathcal{O}$-CM elliptic curve.  
\begin{prop}
there exists an $F$-rational isogeny $\iota: E \ra E'$, where $E'_{/F}$ is an elliptic curve with 
$\mathcal{O}_K$-CM.  Moreover $\iota$ is cyclic of degree $f$.  
\end{prop}
\noindent
This is ``well known'', but lacking a convenient reference we shall sketch the proof.  Over the complex numbers we may 
view $E$ as $\C/\mathcal{O}$, and then the map is just the natural map 
$\C/\mathcal{O} \ra \C/\mathcal{O}_K$.  The rationality of the map over $F$ follows easily from the fact that 
$\mathcal{O}$ is the unique subring of $\mathcal{O}_K$ of index $f$.  
\\ \\
The isogeny $\iota$ induces a homomorphism of Mordell-Weil groups $\iota(F): E(F) \ra E'(F)$.  According to the Proposition, 
the kernel of $\iota(F)$ is $f$-torsion.  Moreover, using the existence of a dual isogeny $\iota^{\vee}: E' \ra E$ 
such that $\iota^{\vee} \circ \iota = [f]$, $\iota \circ \iota^{\vee} = [f]$, one sees that also the cokernel of 
$\iota(F)$ is $f$-torsion.  In particular, if $N$ is an odd prime with $(N,f) = 1$, then 
\[\iota(F): E(F)[N] \stackrel{\sim}{\ra} E'(F)[N]. \]
In particular, if $E$ has an $F$-rational torsion point of order $N$, so does $E'$.  From this it follows that -- 
still for $N$ prime to $f$ -- the least degree of an $\OO(f^2D_0)$-CM point on $X_1(N)$ is at least as large 
as that of an $\OO(D_0)$-CM point on $X_1(N)$.  That is, we have succeeded in reducing Case 2 to Case 1.
\\ \\
\textbf{Case 3: General Case}: Finally suppose we have $D = f^2D_0$ with $N \ | \ f$, and consider an $\OO(D)$-CM elliptic curve $E$ defined over a number field $F$, 
with an $F$-rational $N$-torsion point.  To simplify the analysis, we assume $F$ contains the CM-field $K$ 
(this extra factor of $2$ will not effect the asymptotic analysis).    
\\ \\
The above geometric description of the isogeny $\iota$ shows that $\dim_{\F_N} \ker(\iota) \cap E[N] = 1$, i.e., there exists 
a single point $P_0 \in E[N](\C)$ such that $\langle P_0 \rangle = \ker(\iota) \cap E[N]$.  Consider first any $N$-torsion 
point $P$ which is not in $\langle P_0 \rangle$.  Then $\iota(P)$ is an $F$-rational point on the $\OO(D_0)$-CM elliptic 
curve, i.e., as in Case 2, we immediately reduce to Case 1.  So it suffices to assume that the point $P_0$ is $F$-rational 
and derive lower bounds on $[F:K]$.  
\\ \\
As in $\S 2.3$, Case 3, the mod $N$ Galois representation $\rho_N: \Gal_F \ra GL_2(\Z/N\Z)$ is contained in a 
``pseudo-Cartan subgroup''; taking an ordered basis with $P_0$ as the first vector, we have
\[ \rho(\Gal_F) \subset C_N^{\times} \cong \{ \left[ \begin{array}{cc} a & b \\ 0 & a \end{array} \right] \ 
a \in \F_N^{\times}, \ b \in \F_N \} .\]
So our assumption that $P_0$ is $F$-rational means precisely that 
\[\rho(\Gal_F) \subset \{ \left[ \begin{array}{cc}  1 & b \\ 0 & 1 \end{array} \right] \ b \in \F_N \}. \]
Thus $\det(\rho(\Gal_F)) = 1$, so by Corollary \ref{DETCOR} we deduce $F \supset K(\zeta_N)$.  Now $K(\zeta_N)$ and the 
ring class field $K(j(E))$ are extensions of $K$ of degrees at least $\frac{N-1}{2}$ and $\frac{N-1}{3}$ respectively.  
Moreover, Fact 1e) implies that, loosely speaking, these two extensions are close to being disjoint over $K$, so that 
$K(\zeta_N,j(E))$ has degree at least a universal constant times $(N-1)^2$.  \\ \indent
Let us now see this in more detail: let $E''$ be an elliptic curve 
with $\OO(N^2D_0)$-CM, i.e., with the same CM field but conductor $N$ instead of its multiple $f$.  By class field theory 
$K(j(E)) \subset K(j(E''))$.  But $K(j(E''))$, being the ring class field of conductor $N$, is contained in the $N$-ray 
class field $K(N)$, whereas explicit class field theory shows $\Gal(K(N)/K)$ is a finite abelian group with either $1$ or two generators.  
Therefore the degree of the maximal exponent $2$ abelian subextension of $K(j(E))/K$ is at most $4$.  Combining all estimates, 
we get 
\[ [F:K] \geq [K(j(E),\zeta_N):K] \geq \frac{(N-1)^2}{24}. \]
This is obviously not $O(N)$, so the proof is complete. 

\section{Proof of Theorem \ref{THEOREM7}} \noindent
Here, briefly, is the idea: Start with $E/\Q$ of $j$-invariant $0$.  Enumerate the odd primes $p_n$ which 
are $1$ mod $3$ (hence split in $\Q(\sqrt{-3})$).  Let $K_n$ be the least field over which 
$E$ acquires a point of order $N_n := p_1 \cdots p_n$.  The degree of this field is at most 
\[2 \prod_{i=1}^n (p_i - 1) = 2 \varphi(N_n), \]
and it is known that $\frac{N_n}{\varphi(N_n)} \gg \log \log N_n$.
\\ \\
Proof: Let $K = \Q(\sqrt{-3})$, and $E_{/K}$ an $\OO(-3)$-CM elliptic curve (e.g. $y^2 = x^3 + 1$).  
Let $p_1 < p_2 < \ldots$ be the primes congruent to $1 \pmod 3$, i.e., the primes which split in $K$.  
It follows from the material reviewed in $\S 2.3$ that for each $i$ there is a point $P_i$ on $E$ of order $i$, such that 
$[K(P_i):K] \ | \ (p_i-1)$.  Thus, 
for any positive integer $n$, the field $L_n := K(\{P_i\}_{i=1}^n)$ has a point of order $N_n = p_1 \cdots p_n$ (namely
$P_1 + \ldots + P_n$) and 
\[ d_n := [L_n:K] \leq 2 \prod_{i=1}^n (p_i-1) = 2 \varphi(N_n). \]
Then
\[\frac{|E(L_n)[\tors]|}{d_n}\geq \frac{N_n}{2\varphi(N_n)}, \]
and to complete the proof it is sufficient to establish the following
\\ \\
Claim: There exists $C > 0$ such that for all sufficiently large $n$,
\[\frac{N_n}{2\varphi(N_n)}\geq  C \sqrt{\log(\log(d_n))}.\]
\\ 
The proof of the claim rests on an asymptotic formula due to Mertens, namely\[
\prod_{p\leq x}\frac{1}{1-p^{-1}}\sim e^{-\gamma}\log(x),\]
where the product is taken over all primes less than or equal to $x$, and $\gamma$ is 
Euler's constant \cite[Cor. 6.19]{BD}. From the Prime Number Theorem for Arithmetic Progressions \cite[Thm. 9.12]{BD}, it follows that\[
\prod_{p\leq x,p\equiv  1 (3)}\frac{1}{(1-p^{-1})}\sim e^{-\gamma/2}\sqrt{\log(x)}.\]
Let us now write \[
\frac{N_n}{\varphi(N_n)}=\prod_{i=1}^n\frac{p_i}{p_i-1}=\prod_{p\leq x(n), p\equiv 1 (3)}\frac{1}{1-p^{-1}}.\]
Then we have \[
\frac{N_n}{\varphi(N_n)}\sim e^{-\gamma/2}\sqrt{\log(x(n))}.\]
Again applying the Prime Number Theorem for Arithmetic Progressions, it follows that $\log(x(n))\sim\log(n)$, and also that \[
\log(N_n)=\sum_{i=1}^n\log(p_i)\sim 2\sum_{i=1}^n i\log(i)\sim 2\log{n}
\sum_{i=1}^n i = n(n+1)\log(n).\]
This implies that $\log(\log(N_n))\sim \log(n)\sim \log(x(n))$. Thus \[
\frac{N_n}{\varphi(N_n)}\sim e^{-\gamma/2}\sqrt{\log(\log(N_n))}\geq e^{-\gamma/2}\sqrt{\log(\log(\varphi(N_n)))}
\geq e^{-\gamma/2}\sqrt{\log(\log(d_n/2))},\]
which is sufficient to give the result. 
\\ \\
Remark 8.1: The reader may be wondering whether we could have done better by applying Theorem \ref{JZEROTHM}, which 
says that we can get an $\OO(-3)$-CM point of degree $\frac{p_i-1}{3}$.  However, the factor of $6$ that we gained 
in the proof of this result was via our ability to make a single cyclic twist to get more torsion.  However we cannot 
independently make cyclic twists for each prime $p_i$.  Thus we could improve $d_n$ to $\frac{\varphi(p_1 \cdots p_n)}{3}$ 
but not to $\frac{2}{3^n}\varphi(p_1 \cdots p_n)$.  In fact Serre's Theorem (Theorem \ref{SERRE1}) implies that among 
constructions working with a fixed elliptic curve, or even a fixed $j$-invariant, our lower bound is asymptotically 
optimal.

\section{Proof of Theorem \ref{THEOREM5}}
\noindent
\begin{thm}(Abramovich, \cite{Abramovich})
Let $\Gamma \subset PSL_2(\Z)$ be a congruence subgroup, and $X_{\Gamma} = \Gamma \backslash 
\overline{\mathcal{H}}$ 
the corresponding modular curve.  The gonality of $X_{\Gamma}$ is at least $\frac{7}{800} [PSL_2(\Z):\Gamma].$
\end{thm}
\noindent
Remark 9.1: This result uses results of differential geometry and spectral theory, including an upper 
bound on the leading nontrivial eigenvalue for the Laplacian on the Riemannian manifold $X_{\Gamma}$: 
Abramovich's theorem uses the bound $\lambda_1 \leq \frac{21}{100}$, due to Luo, Rudnick and Sarnak.  
Selberg has conjectured that $\lambda_1 \leq \frac{1}{4}$, which would allow replacement of 
$\frac{7}{800}$ by $\frac{1}{96}$.

\begin{thm}(Faltings, Frey \cite{Frey})
Let $X$ be a curve defined over a number field $K$ with at least one $K$-rational point.  If, for 
any positive integer $d$, $X_{/K}$ has infinitely many points of degree $d$, then 
$\frac{1}{2}\Gon_K(X) \leq d$.
\end{thm}
\noindent
Remark 9.2: The hypothesis is satisfied for all classical modular curves $X_{\Gamma}$ uniformized 
by congruence subgroups of $PSL_2(\Z)$ since such curves always have a cusp rational over their 
``reflex field'' $K$ ($K = \Q$ for the curves $X_1(N)$).  
%For a variant (slightly weaker) result 
%in the case $X(K) = \emptyset$, see \cite{Clark2}.  
\\ \\
When $N$ is prime, the index of $\Gamma_1(N)$ in $PSL_2(\Z)$ is $\frac{N^2-1}{2}$.  Thus we get 
\[\Gon_{\C}(X_1(N)) \geq \frac{7}{1600} (N^2-1) \]
\textbf{unconditionally}, and 
\[\Gon_{\C}(X_1(N)) \geq \frac{1}{192} (N^2-1) \]
\textbf{conditionally} on Selberg's eigenvalue conjecture.
\\ \\%
Therefore we get 
\[\frac{1}{2}\Gon_{\Q}(X_1(N)) \geq \frac{1}{2}\Gon_{\C}(X_1(N)) \geq  \frac{7}{3200}(N^2-1), \]
so if $d \leq [\frac{7}{3200}(N^2-1)] -1$ there are only finitely many points of degree $d$.  
Thus in the statement of Theorem 3 we can take for $C_1$ any constant less than $\frac{7}{3200}$, 
and if Selberg's eigenvalue conjecture holds, we can take any constant less than $\frac{1}{384}$.  
\\ \indent
For part b) we need two facts.  First, for a curve $X$ of genus $g \geq 2$ over any field $K$, one 
can get a degree $2g-2$ map to the projective line by taking an element $f$ of the complete 
linear system associated to the canonical bundle $\Omega^1_{X/K}$, and therefore $\Gon_K(X) \leq 2g(X)-2$. 
Second, for $N > 3$ prime, the genus of $X_1(N)$ is $\frac{N^2-12N+11}{12}$: see e.g. \cite[Theorem 1.1]{JK06}.


\noindent



\begin{thebibliography}{FSWZ90}

\bibitem[Abr96]{Abramovich} D. Abramovich,  \emph{A linear lower bound on the gonality of modular curves}.  
Internat. Math. Res. Notices  1996, no. 20, 1005--1011. 

\bibitem[Ank52]{Ankeny} N.C. Ankeny, \emph{The least quadratic non residue}.  Annals 
of Math. (2) 55 (1952), 65-72.

\bibitem[BD04]{BD} P.T. Bateman and H.G. Diamond,
\emph{Analytic number theory. An introductory course.} World Scientific Publishing Co. Pte. Ltd., Hackensack, NJ, 2004. xiv+360 pp.
%
\bibitem[Bur57]{Burgess} D.A. Burgess, \emph{The distribution of 
quadratic residues and non-residues}.  Mathematika 4 (1957), 106-112.

\bibitem[Cla04]{Clark1} P.L. Clark, \emph{Bounds for torsion on abelian varieties with integral moduli}, $2004$ preprint.

%\bibitem[Cla08]{Clark2} P.L. Clark, \emph{On the Hasse principle for Shimura curves}, Israel Math. J., $2008$.

\bibitem[Coh07]{Cohen} H. Cohen, \emph{Number Theory.  Volume I: Tools and Diophantine Equations}.  Graduate Texts in 
Mathematics 240, Springer-Verlag, 2007. 

\bibitem[Cox89]{Cox} D. Cox, \emph{Primes of the form $x^2+ny^2$.  Fermat, class field theory and complex 
multiplication.}  John Wiley $\&$ Sons, New York, 1989.

%Primes of the form. Fermat, class field theory and complex multiplication.A Wiley-Interscience Publication. John Wiley & Sons, Inc., 
%New York, 1989.

\bibitem[Deu33]{Deuring} M. Deuring, 
\emph{Imagin\"are quadratische Zahlk\"orper mit der Klassenzahl 1.} Math. Z.  37  (1933),  no. 1, 405--415.

\bibitem[Erd52]{Erdos} P. Erd\"os, review of N.C. Ankeny's ``The least quadratic non residue.'' Math reviews $\#$MR0045159.  

\bibitem[Fre77]{Frey} G. Frey, \emph{Some remarks concerning points of finite order on elliptic curves over global fields.}
Ark. Mat. 15 (1977), no. 1, 1--19. 

\bibitem[FSWZ90]{Zimmer2} G. Fung, H. Str\"oher, H. Williams and H. Zimmer, \emph{Torsion groups of 
elliptic curves with integral $j$-invariant over pure cubic fields}. J. Number Theory 36 (1990), 
12-45.

%\bibitem[10]{Gross} B. Gross.  \emph{Arithmetic on Elliptic Curves with Complex Multiplication}.  Lecture 
%Notes in Math. 776, Springer-Verlag, 1980.

\bibitem[HW]{HardyWright} G.H. Hardy, E.M. Wright, \emph{An introduction to the theory of numbers. Fifth edition.} 
The Clarendon Press, Oxford University Press, New York, 1979. xvi+426 pp.

\bibitem[Hei34]{Heilbronn} H. Heilbronn, \emph{On the Class Number in Imaginary Quadratic Fields}. Quart. J. Math. Oxford Ser. 25, 150-160, 1934.

\bibitem[HS99]{HS} M. Hindry and J. Silverman, \emph{Sur le nombre de points de torsion rationnels sur une courbe elliptique.} 
C. R. Acad. Sci. Paris Sér. I Math. 329 (1999), no. 2, 97--100. 

\bibitem[JK06]{JK} D. Jeon and C.H. Kim, \emph{On the arithmetic of certain modular curves}, arXiv preprint, 2006.

\bibitem[JKS04]{JKS} D. Jeon, C. H. Kim and A. Schweizer, \emph{On the torsion of elliptic curves over cubic number
fields}. Acta Arith. 113 (2004) 291–301.

\bibitem[JKP06]{JK06} D. Jeon, C.H. Kim and E. Park, \emph{On the torsion of elliptic curves over quartic number 
fields}. J. London Math. Soc. (2) 74 (2006), 1-12.

\bibitem[Kam86]{Kamienny86} S. Kamienny, \emph{Torsion points on elliptic curves over all quadratic fields}.
Duke Math. J. 53  no. 1 (1986), 157--162. 

\bibitem[Kam92]{Kamienny92} S. Kamienny, \emph{Torsion points on elliptic curves and $q$-coefficients of modular forms}.
Invent. Math. 109 (1992), no. 2, 221--229. 

\bibitem[KM88]{KenkuMomose} M.A. Kenku and F. Momose, \emph{Torsion points on elliptic curves defined over quadratic fields.}
Nagoya Math. J. 109 (1988), 125--149. 

\bibitem[Lan87]{Lang} S. Lang, \emph{Elliptic functions. With an appendix by J. Tate. Second edition.} 
Graduate Texts in Mathematics, 112. Springer-Verlag, New York, 1987.

\bibitem[Maz77]{Mazur} B. Mazur, \emph{Modular elliptic curves and the Eisenstein ideal}, Publ. Math. Inst. Hautes ´ Etudes Sci.
47 (1977) 33–168.
\bibitem[Mer96]{Merel} L. Merel, \emph{Bornes pour la torsion des courbes elliptiques sur les corps de 
nombres}. Invent. Math. 124 (1996), 437-449.

\bibitem[MSZ89]{Zimmer1} H. Miller, H. Str\"oher and H. Zimmer, \emph{Torsion groups of elliptic 
curves with integral j-invariant over quadratic fields}. J. Reine Angew. Math. 397 (1989), 100-161.

\bibitem[Ols74]{Olson} L. Olson, \emph{Points of finite order on elliptic curves with complex 
multiplication}. Manuscripta math. 14 (1974), 195-205.

\bibitem[Par99]{Parent99} P. Parent, \emph{Bornes effectives pour la torsion des courbes elliptiques sur 
les corps de nombres}. 
J. Reine Angew. Math. 506 (1999), 85--116. 

\bibitem[Par03]{Parent} P. Parent, \emph{No 17-torsion on elliptic curves over cubic number fields} .
J. Th\'eor. Nombres Bordeaux 15 (2003), no. 3, 831--838. 

\bibitem[Pari89]{Parish} J.L. Parish, \emph{Rational Torsion in Complex-Multiplication Elliptic Curves}. 
Journal of Number Theory 33 (1989), 257-265.

\bibitem[PWZ97]{Zimmer3} A. Peth\"o, T. Weis and H. Zimmer, \emph{Torsion groups of elliptic curves with 
integral j-invariant over gneeral cubic number fields}.  Internat. J. Algebra Comput. 7 (1997), 353-413.

\bibitem[PY01]{PY} D. Prasad and C.S. Yogananda,
\emph{Bounding the torsion in CM elliptic curves}.  C. R. Math. Acad. Sci. Soc. R. Can.  23  (2001), 1--5.

\bibitem[S66]{Serre66} J.-P. Serre, 
\emph{Groupes de Lie $l$-adiques attach\'es aux courbes elliptiques}. 1966 Les Tendances G\'eom. en Alg\`ebre et Th\'eorie 
des Nombres pp. 239--256 \'Editions du Centre National de la Recherche Scientifique, Paris.

\bibitem[S67]{SerreCM} J.-P. Serre, 
\emph{Complex multiplication}. 1967 Algebraic Number Theory (Proc. Instructional Conf., Brighton, 1965) pp. 292--296 
Thompson, Washington, D.C. 

\bibitem[S72]{Serre72} J.-P. Serre, \emph{Propri\'et\'es galoisiennes des points d'ordre fini des courbes 
elliptiques}.   
Invent. Math.  15  (1972), no. 4, 259--331.

\bibitem[Sie35]{Siegel} C.L. Siegel, \emph{\"Uber die Classenzahl quadtratischer Zahlk\"orper}. Acta Mathematica 1 
(1935), 83-86.

\bibitem[Sbg88]{Silverberg} A. Silverberg, \emph{Torsion points on abelian varieties of CM-type}.  Compositio Math.  68 
(1988),  no. 3, 241--249. 

\bibitem[Sil86]{Silverman} J. Silverman, \emph{The Arithmetic of Elliptic Curves}, Graduate Texts in Mathematics 106, Springer 
Verlag, 1986.

\bibitem[Sil94]{SilvermanII} J. Silverman, \emph{Advanced Topics in the Arithmetic of Elliptic Curves}, Graduate Texts in
Mathematics 151, Springer-Verlag, 1994.

\bibitem[Wat04]{Watkins} M. Watkins, \emph{Class numbers of imaginary quadratic fields.}  
Math. Comp. 73 (2004), no. 246, 907--938. 

\bibitem[Zim76]{Zimmer0} H. Zimmer, \emph{Points of finite order on elliptic curves over number fields}. 
Arch. Math. (Basel) 27 (1976), no. 6, 596--603. 

\end{thebibliography}
\end{document}